\documentclass{amsart}
\usepackage[english]{babel}
\usepackage{graphicx}
\usepackage[hidelinks]{hyperref}
\usepackage{amssymb}
\usepackage{amsmath}
\usepackage[foot]{amsaddr}
\usepackage{amsfonts}
\usepackage{bbm}
\usepackage[dvipsnames, table]{xcolor}
\usepackage[top=1in, bottom=1.25in, left=1.25in, right=1.25in]{geometry}
\usepackage{caption}
\usepackage{subcaption}
\usepackage{multirow}

\usepackage{tikz}
\usetikzlibrary{chains,positioning,shapes.symbols}
\usetikzlibrary{shapes.geometric, arrows}
\usetikzlibrary{shapes,arrows}
\usetikzlibrary{shapes.callouts}
\usetikzlibrary{positioning}
\usepackage{tikzsymbols}

\tikzstyle{startstop} = [rectangle, rounded corners, minimum width=3cm, minimum height=1.2cm, text centered, text width=3cm, draw=black, fill=blue!15]
\tikzstyle{io} = [trapezium, trapezium left angle=70, trapezium right angle=110, minimum width=.8cm, minimum height=1.1cm, text centered,  text width=1.5cm, draw=black, fill=red!40]
\tikzstyle{arrow} = [thick,->,>=stealth]
\tikzstyle{b_arrow_dashed} = [dashed,<->,>=stealth]
\tikzstyle{b_arrow} = [thick,<->,>=stealth]
\tikzstyle{line} = [thick,-,>=stealth]
\tikzstyle{line_dot} = [dotted,-,>=stealth]

\graphicspath{{images/}}

\input{packages.sty}

\definecolor{maria}{HTML}{0090A0}

\colorlet{BLUE}{blue}
\colorlet{BLACK}{black}

\newcommand{\A}[1]{{\color{black}#1}}
\newcommand{\B}[1]{{\color{black}#1}}
\newcommand{\Aa}[1]{{\color{black}#1}}
\newcommand{\Bb}[1]{{\color{black}#1}}
\newcommand{\Aaa}[1]{{\color{black}#1}}

\begin{document}

\title[ROMs for OFCP($\bmu$)]{\Aaa{Space-time POD-Galerkin approach for parametric flow control}}

\author{ Francesco Ballarin$^1$, Gianluigi Rozza$^2$ and  Maria Strazzullo$^{2}$}
\address{$^1$ \textup{Department of Mathematics and Physics, Catholic University of the Sacred Heart, via Musei 41, I-25121 Brescia, Italy.}}
\address{$^2$ \textup{mathLab, Mathematics Area, SISSA, via Bonomea 265, I-34136 Trieste, Italy.}}

\begin{abstract}
In this contribution we propose reduced order methods to fast and reliably solve
parametrized optimal control problems governed by time dependent nonlinear partial differential equations. Our goal is to provide a tool to deal with the time evolution of several nonlinear optimality systems in many-query context, where a system must be analysed for various physical and geometrical features. Optimal control can be used in order to fill the gap between collected data and mathematical model and it is usually related to very time consuming activities:  inverse problems, statistics, etc. Standard discretization techniques may lead to unbearable simulations for real applications. We aim at showing how reduced order
modelling can solve this issue. We rely on a space-time POD-Galerkin reduction in order to solve
the optimal control problem in a low dimensional reduced space in a fast way for several parametric instances. The proposed algorithm is validated with a numerical test based on environmental sciences: a reduced optimal control problem governed by \A{viscous} Shallow Waters Equations parametrized not only in the physics features, but also in the geometrical ones. We will show how the reduced model can be useful in order to recover desired velocity and height profiles more rapidly with respect to the standard simulation, not losing accuracy. \\

\no \textbf{Keywords.} Reduced order modelling, optimal control problems, time dependent nonlinear partial differential equations, Lagrangian approach.
\end{abstract}

\maketitle
\section{Motivations and Historycal Background}

This contribution is rooted in control systems and controllability theory for partial differential equation. A control problem is a system on which you can act through suitable external variables said \emph{controls}  \cite{coron2007control}. The controllability theory answers to the need of steering a system towards \B{a desired configuration}. Is it always possible? Under which conditions can I reach an exact prescribed profile for my system? \\
The problem is quite fascinating and of utmost usefulness in many applications. Even if linear partial differential control systems have many complex aspects and features to be analyzed both theoretically and numerically \cite{coron2007control, GLOWINSKI1992189, lions1968controle,lions1988exact, lions201721}, our main focus will concern nonlinearity in fluid dynamics. In this setting, the problem becomes more challenging and, besides the growing complexity, the need of a control tool increases.  The case of nonlinear partial differential equations is much more complicated to handle.  The control theory for fluid models, e.g. Navier-Stokes equations, prospered in the eighties thanks to the research of J. L. Lions. The main idea of his production relied on the role which nonlinearity plays as a control itself, giving the possibility or preventing the achievement of peculiar motion behaviours \cite{lions144there}.
 This intuition paves the way to a wide range of literature which addresses the problem \cite{agrachev2005navier, coron1996controllability,  lions1969quelques, lions1997controlabilite, lions1998exact, vazquez2008control}.\\
However, in many applied contexts, it is clear that not all the systems are controllable and, furthermore, it is not possible to prove the existence of controls which give the \emph{exact} desired solution profile one wants to reach. This is the reason why the controllabilty theory expands towards optimization and optimal control theory. Namely, the new objective is to find a way to reach \emph{the most similar configuration} with respect to the desired one, satisfying the underling partial differential equation constraint. In the next Section, we will introduce optimal control theory in a parametrized setting, addressing the great importance of a such mathematical model in scientific and engineering contexts.

\section{Introduction}
\label{intro}

This Chapter deals with parametrized optimal control problems (OCP($\boldsymbol{\mu}$)s) governed by time dependent nonlinear parametrized partial differential equations (PDE($\boldsymbol{\mu}$)s). Optimal control is a versatile mathematical tool which has been exploited in many fields of applications: shape optimization, see e.g. \cite{delfour2011shapes,makinen,mohammadi2010applied}, fluid dynamics, see e.g. \cite{dede2007optimal,negri2015reduced,optimal,de2007optimal}, heamodynamics \cite{Ballarin2017,LassilaManzoniQuarteroniRozza2013a,ZakiaMaria,Zakia}, environmental predictions \cite{quarteroni2005numerical,quarteroni2007reduced,Strazzullo1,ZakiaMaria} and more. In a parametrized setting, where a parameter $\boldsymbol{\mu} \in \mathcal  P \subset \mathbb R^d$ could represent physical or geometrical features, \ocp s are exploited in order to study different configurations, which can be  used to better understand the phenomenon one is dealing with.\\
The main motivation for using \ocp s relies in the need of simulations which have to be similar to a given observation or collected data. If on one hand, \ocp s can be useful in many research fields, on the other, they are very complex and demanding from a computational point of view and this issue limits their applicability, most of all in a time dependent framework. Time dependence in \ocp s has been widely described in literature, see e.g. \cite{HinzeStokes,Iapichino2,leugering2014trends,seymen2014distributed,Stoll1,Stoll}. In this context, the required computational resources for \ocp s simulations drastically increase, and parametrized simulations can be unbearable to be performed with standard discretization techniques. To achieve a real-time \ocp $\;$ model, able to reconstruct observable profiles varying with respect to a parameter, a rapid and suitable approximation technique is needed. To this end, we rely on reduced order methods (ROMs). This discretization approach builds a low dimensional framework which can be exploited in order to solve several parametric instances to give real-time information on the model at hand. Working in this reduced space allows us to solve the parametrized optimality system in a small amount of time, by reducing the involved computational costs: for an introduction to the methodology the interested reader may refer to  \cite{antoulas2020interpolatory,benner2017model,hesthaven2015certified,prud2002reliable,RozzaHuynhManzoni2013,RozzaHuynhPatera2008,quarteroni2015reduced}, for example. \\ If we focus on \ocp s applications, there is a wide production concerning steady linear governing equations, see e.g. the following far-from-exhaustive list \cite{bader2016certified,bader2015certified,dede2010reduced,gerner2012certified,Iapichino1,karcher2014certified,karcher2018certified,kunisch2008proper,negri2015reduced,negri2013reduced,quarteroni2007reduced}. Then, the consolidated knowledge about ROM for steady \ocp s has been extended to time dependent \ocp s: moving from \cite{HinzeStokes, Stoll1,Stoll,Yilmaz} as a starting point, in \cite{Iapichino2,Strazzullo2,StrazzulloRB,ZakiaMaria} the main effort is to enlarge the model in a parametrized setting, and generalize the standard algorithm used to build the reduced space framework to time dependency. From now on, we will only focus on POD-Galerkin approach. The motivation relies in its versatility, since it can be even applied to every kind of governing equations: see e.g. \cite{FedeMaria, ZakiaMaria, Zakia} for POD employed in steady \ocp s for nonlinear governing equations, or \cite{Strazzullo3} for the application of such an algorithm to a space-time nonlinear \ocp . In this contribution, the proposed techniques have been tested through parametrized numerical simulations given by distributed control for \A{viscous} Shallow Waters Equations (SWEs), a model capable to simulate coastal current behaviour and used in environmental sciences for monitoring  plans.
\\

We aim at providing a \B{space-time POD-Galerkin strategy for time dependent nonlinear \ocp s}.
This work is outlined as follows. In Section \B{\ref{3}}, \B{we introduce the theoretical formulation for \ocp s following
\cite{hinze2008optimization} and \cite{troltzsch2010optimal}, focusing on time dependent nonlinear problems}. Furthermore, a brief introduction of the space-time algebraic formulation is presented.   Section \B{\ref{4}} deals with the main idea behind reduced order approximation \cite{hesthaven2015certified} and with its \B{application} to space-time nonlinear \ocp s. \B{When dealing with the numerical approximation of time dependent nonlinear \ocp s, we will restrict ourselves to the case of semi-linear PDEs with quadratic nonlinearity in the state variable, to comply with the numerical test presented in Section \ref{5}}: an application in coastal marine management through an \ocp $\,$ governed by SWEs. The proposed numerical test is inspired by \cite{Strazzullo3}, but in this case we consider not only physical, but also geometrical parametrization. Conclusions follow in Section \B{\ref{6}}.

\section{Nonlinear Time Dependent Parametrized Optimal Flow Control Problems}
\label{3}
This Section provides the continuous formulation for nonlinear time dependent \ocp s. We will introduce the Lagrangian approach technique \cite{gunzburger2003perspectives, hinze2008optimization} to minimize a quadratic cost functional constrained to nonlinear time dependent PDE($\bmu$)s. We remark that the analysis is still valid in the Banach spaces setting, but for the sake of clarity, we will restrict ourselves to the simpler case of the \Aa{real} Hilbert spaces.
\subsection{Problem Formulation}
\label{general_problem}
Let us suppose to have a spatial domain $\Omega$ where a physical event described by a time dependent nonlinear PDE($\bmu$) occurs. In order to mathematically represent such an evolution in the time interval $[0, T]$, we define the spaces
\begin{equation}
\mathbb Y = \{y \in L^2(0,T; Y) \text{ such that } y_t \in  L^2(0,T; \B{Y^*}) \} \quad \text{and} \quad
\mathbb Q = L^2(0,T; Y),
\end{equation}
for $Y$ Hilbert space. Indeed, we are provided by a nonlinear \emph{state equation} $G: \mathbb Y \times \mathcal P \rightarrow \mathbb Q\dual$ of the form
\begin{equation}
\label{eq:state}
G(y; \bmu) = f,
\end{equation}
where $y \eqdot y(\bmu) \in \mathbb Y$ is the \emph{state variable}, namely the physical quantity we are interested in, while $f \in \mathbb Q \dual$ is a forcing term and
$\mathcal P \subset \mathbb R^P$
is a parameter space of dimension $P \geq 1$, which can describe physical and/or geometrical features of the system at hand. Moreover, we call $\mathcal L(\cdot, \cdot)$ the space of the continuous \B{linear} functions between two spaces. The considered PDE($\bmu$)s present the following form: \begin{equation}
\label{eq:state_form}
G(y; \bmu) =
y_t + E_{\textit{n}\ell}(y; \bmu) + E_{\ell}(y; \bmu)
,
\end{equation}
where $E_{\ell} \in \Cal L(\state, \mathbb Q \dual)$ and $E_{\textit{n}\ell}$ represent the linear and nonlinear contributions to the equation, respectively. The term $y_t$ describes the time evolution of the equation. \\ We now want to change the behaviour of the state variable steering it to a desired profile, say $y_\text{d} \eqdot y_\text{d}(\bmu) \in \mathbb Y_{\text{obs}} \supseteq \mathbb Y$, thanks to the action of a \emph{control variable} $u \eqdot u(\bmu) \in \mathbb U$, where $\mathbb U = L^2(0,T; U)$ with $U$ another Hilbert space. Thanks to all the previous definitions, we are able to define the \emph{controlled equation} $\mathcal E(y,u; \bmu)$, where
$\mathcal E: \mathbb Y \times \mathbb U \times \mathcal P \rightarrow \mathbb Q\dual$ of the form:
\begin{equation}
\label{eq:control_eq}
\mathcal E(y,u; \bmu) \eqdot\; G(y; \bmu) - C(u) - f = 0,
\end{equation}
where $C \in \Cal L(\control, \mathbb Q \dual)$ is related to the control variable action on the system in order to change the original state variable.
The \ocp $\;$ reads: given a $\bmu \in \mathcal P$, find the pair $(y,u) \in \mathbb Y \times \mathbb U$ which solves
\begin{equation}
\label{min_problem}
\min_{y \in \mathbb Y_{\text{ad}} \subset \mathbb Y, u \in \mathbb U_{\text{ad}} \subset \mathbb U} J(y,u; y_\text{d}) \text{ subject to } \mathcal E(y,u; \bmu) = 0,
\end{equation}
where $J: \mathbb Y \times \mathbb U \times \mathbb Y_{\text{obs}} \rightarrow \mathbb R$ is the \emph{cost functional} defined by
\begin{equation}
J(y,u; y_\text{d}) \eqdot \half \norm{y - y_\text{d}}_\mathbb {Y_{\text{obs}}}^2 + \alf \norm{u}_{\mathbb U}^2,
\end{equation}
and $\alpha \in (0, 1]$ is a \emph{penalization parameter}. It is clear that the smaller is the value of $\alpha$, the more the control variable will influence the system. Problem \eqref{min_problem} admits a solution if \cite[Section 1.5.2]{hinze2008optimization}:
\begin{enumerate}[(i)]
\item $\mathbb U_{\text{ad}}$ is convex, bounded and closed;
\item $\mathbb Y_{\text{ad}}$ is convex and closed;
\item for every $\bmu \in \Cal P$, the controlled system $\mathcal E (y,u; \bmu) = 0$ has a bounded solution map
$u \in \mathbb U \mapsto y(u) \in \mathbb Y$;
\item for a given $\bmu \in \Cal P$, the map $(y,u, \bmu) \in \mathbb Y \times \mathbb U \times \mathcal P \rightarrow \mathcal E (y,u; \bmu) \in \mathbb Q\dual$ is weakly continuous with respect to (w.r.t.) the first two arguments;
\item for a given $y_{\text{d}} \in \mathbb Y_{\text{obs}}$, the \B{cost} functional $J(y,u; y_\text{d})$ is weakly lower semicontinous w.r.t.\ the first two arguments.
\end{enumerate}
To apply the Lagrangian theory, we define $z \eqdot z(\bmu) \in {\mathbb Y^{\ast \ast}} = \mathbb Y \subset \mathbb Q$ be an arbitrary \emph{adjoint variable}, which allows to translates the \eqref{min_problem} in an unconstrained minimization problem. Indeed, calling with $X \eqdot X(\bmu) = (y(\bmu),u(\bmu),z(\bmu)) \in \mathbb X \eqdot \mathbb Y \times \mathbb U \times \mathbb Y$ the global variable of the system at hand, we can build the \emph{Lagrangian functional}
$\Lg:\mathbb X \times \mathbb Y_{\text{obs}} \times \mathcal P \rightarrow \mathbb R$ as
\begin{equation}
\label{lagrangian_functional}
\Lg(X; y_\text{d}, \bmu) \eqdot J(y, u; y_\text{d}) + \la z, \Cal E(y,u; \bmu) \ra_{\mathbb Q \mathbb Q\dual},
\end{equation}
where with $\la \cdot, \cdot \ra_{\mathbb Q\mathbb Q\dual}$ we mean the duality pairing between $\mathbb Q$ and $\mathbb Q\dual$.  Let us assume that the following hold:
\begin{enumerate}[resume*]
\item $\mathbb U$ is nonempty;
\item $J : \mathbb Y \times \mathbb U \times \mathbb Y_{\text{obs}} \rightarrow \mathbb R$ and $\mathcal E : \mathbb Y \times \mathbb U \times \mathcal P \rightarrow \mathbb Q\dual$ are continuously Fr\'echet differentiable w.r.t.\ the first two arguments;
\item given $\bmu \in \mathcal P$, the controlled system $ \mathcal E(y, u; \bmu) = 0$ has a unique solution $y = y(u) \in \mathbb Y$ for all $u \in \mathbb U$;
\item given $\bmu \in \mathcal P$, $D_y \mathcal E (y, u; \bmu) \in \mathcal L(\mathbb Y, \mathbb Q\dual)$ has a bounded inverse for all control variables $u$.
\end{enumerate}
From (ix), it is clear that $D_{\star}$ indicates the Fr\'echet derivative w.r.t.  a variable $\star$ and we will exploit this notation from now on. Thanks to hypotheses
(vi) - (ix), given a solution $({y}, u) \in \state \times \control$ of \eqref{min_problem} for a given $\bmu \in \Cal P$, there exists an adjoint variable $z \in \state$ which satisfies the following \emph{optimality system} \cite{hinze2008optimization}:
\begin{equation}
\label{KKT}
\begin{cases}
D_{y}\Lg(X; y_\text{d}, \bmu) [\omega] = 0 & \forall \omega \in \mathbb Q,\\
D_u\Lg(X; y_\text{d}, \bmu) [\kappa] = 0 & \forall \kappa \in \control,\\
D_z\Lg(X; y_\text{d}, \bmu) [\zeta] = 0 & \forall \zeta \in \mathbb Q,\\
\end{cases}
\end{equation}
or equivalently, in strong form
\begin{equation}
\label{KKT_explicit}
\begin{cases}
y + D_y \mathcal E (y, u; \bmu)\dual (z) = y_\text{d}, &\\
\alpha u - C\dual (z )= 0, &\\
\mathcal E(y,u;\bmu) = 0, &\\
\end{cases}
\end{equation}
where $D_y \mathcal E (y, u; \bmu) \dual \in \Cal L(\mathbb Q, \state \dual)$ is the adjoint operator of the Fr\'echet linearization of $\mathcal E (y,u; \bmu)$ w.r.t the state variable, while $C\dual \in \Cal L(\mathbb Y, \mathbb U\dual) $ is the adjoint of the control operator.
Furthermore, the system \eqref{KKT} can be recast in compact form: given $\bmu \in \Cal P$, find $X \in \mathbb X$ such that
\begin{equation}
\label{ocp}
\mathcal G(X; \bmu) = {\Cal F},
\end{equation}
with
\begin{equation*}
\mathcal G(X; \bmu) \eqdot \begin{bmatrix} y + D_y \mathcal E (y, u; \bmu)\dual (z) \\ \alpha u - C\dual (z) \\ G(y, \bmu) - C(u) \end{bmatrix} \quad
\text{and} \quad \Cal F \eqdot \begin{bmatrix} y_{\text{d}} \\ 0 \\ f \end{bmatrix}.
\end{equation*}
We underline that the dual variable $z$ is considered in $\mathbb Y$ in order to guarantee a proper definition of the optimality system due to the presence of a backward time evolution of the form $-z_t$ in the expression of $D_y \mathcal E (y, u; \bmu) \dual$. \\
In this context, we will always assume that the state equation is local invertible for every parametric instance, i.e. assumptions (viii) and (ix) always hold. We underline that in the nonlinear case it may happen to find multiple solutions for a given parameter. However, it will not be the case of this work, since we are  restricting ourselves to a well-posed setting for the state equation. \B{Indeed, we stress that in the nonlinear case it is only possible to recover existence results for the system \eqref{min_problem}. Uniqueness results are strongly related to the state equation one is dealing with and multiple solutions for the optimality system \eqref{ocp} can be found, see for example \cite{FedeMaria}}.

In the next Section, we will show the space-time approximation of the optimality system at hand, generalizing the strategy already presented in \cite{Glas2017,urban2012new,yano2014space,yano2014space1}.  \B{We will focus on the case of semi-linear governing equations with quadratic dependence in the state variable, guided by the numerical results provided in Section \ref{5}}.

\subsection{The Space-Time Approximation}
\label{FE}
The next step to be taken is the numerical approximation of the optimal solution $X$ of \eqref{ocp} in order to investigate its features varying $\bmu \in \mathcal P$ in a space-time fashion. This approximation approach is a versatile tool already exploited in the discretization of parabolic equations and several \ocp s, see e.g. \cite{Glas2017, HinzeStokes, hinze2008optimization, HinzeNS, Langer2020, Strazzullo2, StrazzulloRB, Strazzullo3, yano2014space, yano2014space1, urban2012new}. We first focus on the space discretization, which is performed through the Finite Element (FE) technique.  Let us define the FE function space $Y^{N_{\text{FE}}^y} = Y \cap\mathbb K_{r_y}$ and $U^{N_{\text{FE}}^u} =  U \cap\mathbb K_{r_u}$, where
$$
\mathbb K_r = \{ v \in C^0(\overline \Omega) \; : \; v |_{K} \in \mbb P^r, \; \; \forall \, K \in \Cal{T} \},
$$
with $\mbb P^r$ is the space of all the polynomials of degree at most equal to $r$ and  $K$ is an element of a triangularization $\Cal T$ of the spatial domain $\Omega$.  We can now consider the {semi-discrete} function spaces
 $ \mathbb Y^{N_{\text{FE}}^y} = \Big \{
y \in L^2(0,T; Y^{N_{\text{FE}}^y}) \;\; \text{s.t.} \;\;  \dt{y} \in L^2(0,T; {(Y^{N_{\text{FE}}^y})}^{\ast})\Big \}$, $\mathbb Q^{N_{\text{FE}}^y} = L^2(0,T; Y^{N_{\text{FE}}^y})$ and $\mathbb U^{N_{\text{FE}}^u} = L^2(0,T; U^{N_{\text{FE}}})$.
 Once made the spatial discretization step, a time discretization over the interval $[0,T]$ must be taken into account resulting in the final space-time discrete spaces, which are denoted by $\mathbb Y^{N_{\text{FE}}^y}_{N_t}$, $\mathbb Q^{N_{\text{FE}}^y}_{N_t}$ and $\mathbb U^{N_{\text{FE}}^u}_{N_t}$, where $N_t$ is the number of the considered timesteps. \B{Following \cite{Langer2020}, we decided to use the same space-time approximation for $\mathbb Y^{N_{\text{FE}}^y}_{N_t}$ and $ \mathbb Q^{N_{\text{FE}}^y}_{N_t}$}. For the sake of notation, we will refer to the space-time function spaces as $\mathbb Y \discy \equiv \mathbb Q \discy$ and $\mathbb U \discu$ and, as a consequence, $\mathbb X\disc \eqdot \mathbb Q \discy \times \mathbb U \discu \times \mathbb Q \discy$, i.e. $\mathcal N = 2\mathcal N_y + \mathcal N_u $, with $\mathcal N_y = N_{\text{FE}}^y \cdot N_t$ and $\mathcal N_u = N_{\text{FE}}^u \cdot N_t$. In this finite dimensional setting the problem to be solved reads: given $\bmu \in \Cal P$ and observation $y_{\text{d}}\disc\in \mathbb Q \discy$ find $X\disc \eqdot X \disc (\bmu) \in \mathbb X\disc$ such that
\begin{equation}
\label{FE_ocp}
\mathcal G(X\disc; \bmu) = {\mathcal F}.
\end{equation}

In the following, we make clear the algebraic structure we exploited in the numerical experiment presented in Section \ref{SWE}. First of all, let us divide the time interval $[0,T]$ in $N_t$ equispaced subintervals of length $\Delta t$ and let us call $t_k = k\Delta t$ for $k = 0, \dots, N_t$ a generic time instance. The variables $y\disc_k$m $u\disc_k$ and $z\disc_k$ for a specific timestep, can be represented with FE basis $\{\phi^i\}_{i=1}^{N_{\text{FE}}^y}$ and $\{\psi^i\}_{i=1}^{N_{\text{FE}}^u}$ for $Y^{N_{\text{FE}}^y}$ and $U^{N_{\text{FE}}^u}$, respectively as follows
$$
y\disc_k = \sum_{i=1}^{N^y_{\text{FE}}}y^{i}_k \phi^{i}, \hspace{1cm}
u\disc_k = \sum_{i=1}^{N^u_{\text{FE}}}u^{i}_k \psi^{i}, \hspace{1cm} \text{and} \hspace{1cm}
z\disc_k = \sum_{i=1}^{N^y_{\text{FE}}}z^{i}_k \phi^{i}.
$$

We now define the space-time state, control and adjoint vectors  $$\mathsf y = 
\begin{bmatrix}
\bar y_1 \\
\vdots\\
 \bar y_{N_t}
\end{bmatrix}, \quad 
\mathsf u = 
\begin{bmatrix}
\bar u_1 \\
\vdots\\
 \bar u_{N_t}
\end{bmatrix}, 
\quad \text{and}  \quad
\mathsf z = 
\begin{bmatrix}
\bar z_1 \\
\vdots\\
 \bar z_{N_t}
\end{bmatrix}.$$ 
\B{Namely, $\bar y_k$, $\bar u_k$ and $\bar z_k$ are the column vectors with FE coefficients of the variables at time instance $t_k$, with $k=1, \dots, N_{t}$, i.e. 
$$\bar y_k = 
\begin{bmatrix}
y_k^1 \\
\vdots\\
 y_k^{N_{FE}^y}
\end{bmatrix}, \quad 
\bar u_k = 
\begin{bmatrix}
u_k^{1} \\
\vdots\\
u_k^{N_{FE}^u}
\end{bmatrix}, 
\quad \text{and}  \quad
\bar z_k = 
\begin{bmatrix}
z_k^1 \\
\vdots\\
 z_k^{N_{FE}^y}
\end{bmatrix}.$$
Applying the same strategy to the initial time condition, to the desired state and the forcing term, we can define $$\mathsf y_0 = 
\begin{bmatrix}
\bar y_0 \\
0 \\
\vdots\\
 0
\end{bmatrix}, 
\quad
\mathsf y_d = 
\begin{bmatrix}
\bar {y_d}_1\\
\bar {y_d}_2\\
\vdots\\
\bar {y_d}_{N_t}
\end{bmatrix}, 
\quad \text{and} \quad 
\mathsf f =
\begin{bmatrix}
\bar {f}_1\\
\bar {f}_2\\
\vdots\\
\bar {f}_{N_t}
\end{bmatrix}, 
$$ respectively, where $\bar y_0$, $\bar y{_{d}}_{k}$ and $\bar f_{k}$, in analogy with the aforementioned space-time variables, are the column vectors given by the FE coefficients in their respectively function spaces,
with  $k=1, \dots, N_{t}$.} 
We begin the discretization analysis from the state equation, for the sake of clarity. Concerning the space discretization level, once applied the controlled equation to the FE basis, we can derive the matrices
$\Aa{\mathsf {E}_{\textit{n}\ell}(\mathsf y; \bmu)} + \mathsf {E}_{\ell}(\bmu) - \mathsf C$, omitting for the moment the time evolution. Furthermore, we need the mass matrices $\mathsf M_y$ and $\mathsf M_u$ for state/adjoint variables and control, respectively.  At every time instance, after performing a backward Euler discretization in time, we have to solve the following system:
\begin{equation}
{\mathsf M_y} \bar y_k + \Delta t \underbrace{(\Aa{\mathsf {E}_{\textit{n}\ell}(\mathsf y; \bmu)}  + \mathsf {E}_{\ell}(\bmu) )}_{\Aa{\mathsf E(\mathsf y, \bmu)}}\bar y_k - {\Delta t}{\mathsf C} \bar u_k= {\mathsf M_y} \bar y_{k-1}  + \bar f_k \Delta t
\spazio \text{for } 1 \leq k \leq N_t.
\end{equation}
{Thus}, the whole space-time system reads
\begin{equation}
\label{test}
\underbrace{
\begin{bmatrix}
{\mathsf M_y} + \Delta t \Aa{\mathsf E(\mathsf y, \bmu)} &  & & \\
- {\mathsf M_y} & {\mathsf M_y} +\Delta t \Aa{\mathsf E(\mathsf y, \bmu)} &  & \\
& - {\mathsf M_y} & {\mathsf M_y}+ \Delta t \Aa{\mathsf E(\mathsf y, \bmu)} &  & \\
&  & \ddots & \ddots & \\
&  &  & -{\mathsf M_y} & {\mathsf M_y} + \Delta t \Aa{\mathsf E(\mathsf y, \bmu)}  \\
\end{bmatrix}
}_{\Aa{\mathsf K(\mathsf y; \bmu)}}
\begin{bmatrix}
\bar y_1\\
\bar y_2\\
\bar y_3 \\
\vdots \\
\bar y_{N_t}
\end{bmatrix}
$$
$$
\qquad \qquad \qquad
\qquad \qquad \qquad
- {\Delta t}\begin{bmatrix}
\mathsf C&  & & \\
&\mathsf C   &  & \\
&  & \mathsf C   &  & \\
&  &  &  \ddots & \\
&  &  &  & \mathsf C   \\
\end{bmatrix}
\begin{bmatrix}
\bar u_1\\
\bar u_2\\
\bar u_3 \\
\vdots \\
\bar u_{N_t}
\end{bmatrix}
\hspace{-1mm}= \hspace{-1mm}
\begin{bmatrix}
{\mathsf M_y}\bar y_0 + \Delta t \bar f_1\\
0 + \Delta t  \bar f_2\\
0 + \Delta t  \bar f_3\\
\vdots \\
0 + \Delta t  \bar f_{N_t}
\end{bmatrix}.
\end{equation}
{We underline that the system presents a nonlinear dependence from the state in the term $\mathsf{E}_{n\ell}(\bmu)$, related to the nonlinear contribution of the state operator $E_{n\ell}(y;\bmu)$. }
The space-time state equation can be written in compact form as
\begin{equation}
\Aa{\mathsf K(\mathsf y; \bmu)} \mathsf  y - {\Delta t} \mathsf C_{\text{st}} \mathsf u = \mathsf {M_y}_{\text{st}} \mathsf y_0 + \Delta t  \mathsf f,
\end{equation}
where the subscript  ``st'' indicates the all-at-once matrices, namely $\mathsf C_{\text{st}}$ is the block-diagonal matrix which entries are given by $\mathsf C$ of dimension  $ \mathbb R^{\mathcal N_y} \times \mathbb R^{\mathcal N_u}$ and $\mathsf {M_y}_{\text{st}}$ is the matrix made by $\mathsf M_y$ on the diagonal, of dimension  $ \mathbb R^{\mathcal N_y} \times \mathbb R^{\mathcal N_y}$. \\
We now take into account the optimality equation, which easily reads:
\begin{equation}
\B{\alpha \Delta t \mathsf M_u \bar u_k - \Delta t \mathsf C^T \bar z_k = 0} \spazio \text{for } 1 \leq k \leq N_t,
\end{equation}
or in compact form
\begin{equation}
\B{\alpha \Delta t \mathsf {M_u}_{\text{st}} \mathsf u-  \Delta t \mathsf C^T_{\text{st}} \mathsf z = 0},
\end{equation}
where $\mathsf {M_u}_{\text{st}} $ is the block-diagonal matrix which entries are given by $\mathsf M_u$ of dimension  $ \mathbb R^{\mathcal N_u} \times \mathbb R^{\mathcal N_u}$.
 We can now take into account the adjoint equation. To this purpose, we have to explicit the algebraic structure of $D_y \Cal E(y,u; \bmu)$. The assumption of the quadratic nonlinearity  in the state variable gives the following form to the Fr\'echet derivative of the controlled state equation with respect to the state $y$:
$\mathsf {E}_{\textit{n}\ell}'[\mathsf y] (\bmu)+ \mathsf E_{\ell}(\bmu)$,  where the linear state structure remains the same, while the nonlinear operator is linearized in $\mathsf {E}_{\textit{n}\ell}'[\mathsf y] $. It is clear that the control operator $\mathsf C$ has to disappear since it does not depend on the state variable. Thus, performing a forward Euler method which is equivalent to an implicit scheme due the backward parabolic nature of the adjoint equation, at each time instance it reads:
\begin{equation}
{\mathsf M_y} \bar z_{k-1}  = {\mathsf M_y} \bar z_k  + \Delta t ( - {\mathsf M_{\text{obs}}} \bar y_{k-1}\underbrace{ - \mathsf {E}_{\textit{n}\ell}'[\mathsf y]^T(\bmu) - \mathsf E_{\ell}^T(\bmu)}_{-{\mathsf E^{\text{adj}}}^T(\bmu)} \bar p_{k-1} + {\mathsf M_{\text{obs}}} \bar y_{d_{k-1}} )
\spazio \text{for } 1< k \leq N_t,
\end{equation}
where $\mathsf M_{\text{obs}}$ is the state mass matrix restricted to the observation domain. {We stress that $\mathsf E^{\text{adj}}(\bmu)$ hides the state variable dependence. }Futhermore, also in this case, we can write the whole all-at-once system:
\begin{equation*}
\underbrace{
\begin{bmatrix}
{\mathsf M_y}+ \Delta t {\mathsf E^{\text{adj}}}^T(\bmu) & - {\mathsf M_y}& & \\
&{\mathsf M_y} + \Delta t {\mathsf E^{\text{adj}}}^T(\bmu) & -{\mathsf M_y}  & \\
&  & \ddots & \ddots  & \\
&  &  & {\mathsf M_y} + \Delta t {\mathsf E^{\text{adj}}}^T(\bmu) & -{\mathsf M_y}  \\
&  &  &  & {\mathsf M_y} + \Delta t {\mathsf E^{\text{adj}}}^T(\bmu) \\
\end{bmatrix}
}_{{\mathsf K}'(\bmu)}
\begin{bmatrix}
\bar z_1\\
\bar z_2\\
\bar z_3 \\
\vdots \\
\bar z_{N_t}
\end{bmatrix}
$$
$$
+
\Delta t
\begin{bmatrix}
{\mathsf M_{\text{obs}}}&  & & \\
& {\mathsf M_{\text{obs}}}&  & \\
&  & {\mathsf M_{\text{obs}}}   &  & \\
&  &  & \ddots & \\
&  &  &  & {\mathsf M_{\text{obs}}}   \\
\end{bmatrix}
\begin{bmatrix}
\bar y_1\\
\bar y_2\\
\bar y_3 \\
\vdots \\
\bar y_{N_t}
\end{bmatrix}
=
\begin{bmatrix}
\Delta t {\mathsf M_{\text{obs}}}\bar y_{d_1} \\
\Delta t {\mathsf M_{\text{obs}}}\bar y_{d_2} \\
\Delta t {\mathsf M_{\text{obs}}}\bar y_{d_3} \\
\vdots \\
\Delta t {\mathsf M_{\text{obs}}} \bar y_{d_{N_t}}
\end{bmatrix}.
\end{equation*}
Then, in a more compact notation, the adjoint equation reads:
\begin{equation}
\mathsf K'(\bmu)\mathsf z + \Delta t {\mathsf M_{\text{obs}}}_{\text{st}} \mathsf y = \Delta t {\mathsf M_{\text{obs}}}_{\text{st}} \mathsf y_{\text{d}},
\end{equation}
where { ${\mathsf M_{\text{obs}}}_{\text{st}} \in \mathbb R^{\mathcal N_y} \times \mathbb R^{\mathcal N_y}$ is the block-diagonal matrix which entries are given by $\mathsf M_{\text{obs}}$.}
Finally, we are able to collect all the information and build the global all-at-once system:
\begin{equation}
\label{one_shot}
\overbrace{
\begin{bmatrix}
\Delta t {\mathsf M_{\text{obs}}}_{\text{st}}&  0 &\mathsf {K}'(\bmu) \\
0 & \alpha \Delta t {\mathsf M_{\mathsf u}}_{\text{st}}&  -\Delta t {\mathsf C}_{\text{st}}^T \\
\Aa{\mathsf K(\mathsf y; \bmu)} & -{\Delta t} {\mathsf C}_{\text{st}} & 0 \\
\end{bmatrix}
\underbrace{
\begin{bmatrix}
\mathsf y \\
\mathsf u \\
\mathsf z \\
\end{bmatrix}}
_{\mathsf X}}^{\mathsf G(\mathsf X; \bmu)}
=
\overbrace{
\begin{bmatrix}
\Delta t {\mathsf M_{\text{obs}}}_{\text{st}} \mathsf y_d \\
0 \\
{\mathsf {M_{{\mathsf y}_{\text{st}}}}} \mathsf y_0 + \Delta t \mathsf f \\
\end{bmatrix}}^{\mathsf F}.
\end{equation}
{We remark that the nonlinearity of $\mathsf G(\mathsf X; \bmu)$ derives from $\Aa{\mathsf K(\mathsf y; \bmu)}$ and $\mathsf K'(\bmu)$, due to the nonlinear terms $\Aa{\mathsf E_{n\ell}(\mathsf y; \bmu)}$ and $\mathsf E_{n\ell}'[\mathsf y]^T(\bmu)$. However, for the sake of notation, we omitted the direct $\mathsf X-$dependence from the matrices.}\\ 
The nonlinear system, then, can be recast in residual formulation as
\begin{equation}
\label{G_compact_FE}
\mathsf R(\mathsf X; \bmu)\eqdot \mathsf G(\mathsf X; \boldsymbol \mu) - \mathsf F = 0,
\end{equation}
where $\mathsf R(\mathsf X; \bmu)$ will be called the \emph{global residual} of the optimality system.
To solve system \eqref{G_compact_FE}, we employed Netwon's method: namely we iteratively solve
\begin{equation}
\mathsf {X}^{j + 1} = \mathsf {X}^j+ \mathsf{Jac}(\mathsf X^{j}; \boldsymbol \mu)^{-1}(\mathsf F - \mathsf G(\mathsf X^j; \boldsymbol \mu)), \spazio j \in \mathbb N,
\end{equation}
until a residual based convergence criterion is satisfied.
We recall that the matrix $\mathsf {K}'(\bmu)$ still depends on the state vector $\mathsf y$ in the term $\mathsf E_{\textit{n}\ell}'[\mathsf y]^T$. Then, the linearization w.r.t. $\mathsf y$ of $\mathsf K'(\bmu)$ leads to a new term in the formulation:
\begin{equation}
\mathsf D_{\mathsf y}( \mathsf {E}_{\textit{n}\ell}'[\mathsf y^j]^T)[\mathsf z^j].
\end{equation}
Namely, the Jacobian matrix will have the following form:
\begin{equation}
\label{J_ocp}
\mathsf{Jac}(\mathsf X^j; \bmu) = \begin{bmatrix}
\Delta t \mathsf {M_{\text{obs}}}_{\text{st}} + \mathsf D_{\mathsf y}( \mathsf {E}_{\textit{n}\ell}'[\mathsf y^j]^T)[\mathsf z^j] & 0 & \mathsf K^{\text{l}}(\Aa{\mathsf y_j}; \bmu)^T\\
0 & \alpha \Delta t {\mathsf M_u}_{\text{st}} & - \Delta t{\mathsf C}_{\text{st}}^T \\
\mathsf K^{\text{l}}(\Aa{\mathsf y_j}; \bmu)& - \Delta t \mathsf {\mathsf C}_{\text{st}} & 0 \\
\end{bmatrix},
\end{equation}
where each matrix taken into consideration is now linear in the $j-$th value of one of the involved variables and $\mathsf K^{\text{l}}(\Aa{\mathsf y_j};\bmu)$ is the linearized version of $\Aa{\mathsf K(\mathsf y,\bmu)}$. Thanks to this remark we are able to show  the \emph{saddle point} structure of the system at hand. Indeed, equation \eqref{J_ocp} can be written as \begin{equation}
\label{J_saddle}
\mathsf{Jac}(\mathsf {X}^j; \bmu) =
\begin{bmatrix}
\mathsf A & \mathsf B^T \\
\mathsf B & 0 \\
\end{bmatrix},
\end{equation}
where
\begin{equation}
\mathsf A =
\begin{bmatrix}
\Delta t \mathsf {M_{\text{obs}}}_{\text{st}} + \mathsf D_{\mathsf y}( \mathsf {E}_{\textit{n}\ell}'[\mathsf y^j]^T)[\mathsf z^j]  & 0\\
0 & \alpha \Delta t {\mathsf M_u}_{\text{st}}  \\
\end{bmatrix}
\text{ and } \spazio \mathsf B =
\begin{bmatrix}
\mathsf {K}^{\text{l}} (\Aa{\mathsf y_j};\bmu) & - \Delta t \mathsf C_{\text{st}}
\end{bmatrix}.
\end{equation}
The proposed framework is very common in many CFD applications, from Stokes equations to PDEs contrained optimization. \B{In the saddle point setting, following the Brezzi theory \cite{brezzi1974existence}, to guarantee the existence of a unique solution of \eqref{J_saddle}, the matrix $\mathsf A$ should be invertible and the following \emph{Brezzi inf-sup condition} should be verified \cite{benzi_golub_liesen_2005,brezzi1974existence}}:
\begin{equation}
\label{FE_lbb}
\beta \disc(\bmu) \eqdot \adjustlimits \inf_{0 \neq \mathsf z} \sup_{0 \neq \mathsf x} \frac{\mathsf z^T\mathsf B \mathsf x}{\norm{\mathsf x}_{\state \times \control}\norm{\mathsf z}_{ \state}} \geq \hat{\beta}(\bmu)\disc > 0,
\end{equation}
where $\mathsf x =[ 
\mathsf y,
\mathsf u
]^T
$. \B{For linear state equations hypotheses (i)-(ix) assure the well-posedness of the saddle point structure \cite{hinze2008optimization, negri2015reduced, negri2013reduced, Strazzullo2}. This is not the case for nonlinear settings, where the where the fulfillment of  (i)-(ix) does not ensure uniqueness and the Brezzi theorem should be verified from case to case.} In the FE context, the inequality \eqref{FE_lbb} holds when the function spaces for state and adjoint coincide \cite{negri2015reduced,negri2013reduced}. The assumption $z \in \state$, will guarantee the fulfillment of the inf-sup stability condition in the space-time approximation once provided at the continuous level.
It is clear that in order to solve the all-at-once optimization problem in a parametrized setting, for a given $\boldsymbol \bmu \in \Cal P$, we have to deal with a high-dimensional systems. The solution of many-query and/or real-time tasks require growing computational resources and computational time for simulations.
\A{In this context, the space-time approximation has some limitations most of all when one relies, as in our case, on the direct solution of the optimality system. To lighten this issue, besides the employment of more computational resources, proper preconditioners and multigrid approaches cen be used, see e.g. \cite{benzi_golub_liesen_2005,schoberl2007symmetric,Stoll1,Stoll} and the references therein.}  In the next Section we will introduce ROMs for space-time nonlinear \ocp s, providing a general approximation strategy which can solve the issue of the huge amount of computational costs that are usually associated to a standard space-time solution process. 

\section{ROMs for Nonlinear Space-Time \ocp s}
\label{4}

This Section introduces ROMs for nonlinear space-time \ocp s. We refer to \cite {Strazzullo1, ZakiaMaria, Zakia} for previous contributions to ROM for nonlinear \ocp s and to \cite{Strazzullo3} for their extension to time dependent nonlinear governing equations. First, we introduce the ROMs ideas and we will briefly focus on the standard approaches to make the strategy effiecient in terms of computational resources. Then, we will describe the POD-Galerkin basis construction algorithm, see \cite{ballarin2015supremizer, benner2017model, burkardt2006pod, Chapelle2013, hesthaven2015certified} as general references. We will exploit the classical aggregated spaces technnique, following the linear \ocp s fashion, as already presented in \cite{bader2016certified,bader2015certified,dede2010reduced,gerner2012certified,karcher2014certified,karcher2018certified, negri2015reduced,negri2013reduced,quarteroni2007reduced}.  \Aaa{We underline that the proposed strategy is strictly related to the linear quadratic case, i.e.\ \ocp s governed by linear PDE($\bmu$)s with a quadratic cost functional. However, it represents a classic choice to deal with more general nonlinear frameworks, see e.g.\ \cite{quarteroni2015reduced}.}

\subsection{General Reduction Strategy}
\label{sec_rom_gen}
In Section \ref{intro}, we introduced the importance of parametric optimal control in several field of applications. Parameters can represent physical features and/or geometrical ones and, in many-query and real-time contexts, there exists the need to study several parametric instances to better understand the properties of a system. For this task,  the space-time formulation can be unbearable and limit the knowledge capability of \ocp s due to the huge amount of computational resources needed for their simulations. The ROM aim at building a \emph{low-dimensional} surrogate function space, in order to decrease the needed time for a simulation, guaranteeing a better parametric analysis in a small amount of time. 
We now provide the main ROM ideas for \ocp s. Let us consider the global variable
$X(\boldsymbol \mu) = (y(\boldsymbol \mu), u(\boldsymbol \mu),z(\boldsymbol \mu))$, parametric solution of \eqref{ocp}. In this Section we make the parameter dependence explicit, for the sake of clarity. Indeed, it will be useful to understand the basics and main features of this discretization approach. 
The first phase of the ROM relies in the construction of basis functions to represent the \emph{high fidelity solution} $\Aa{X^{\Cal N}(\bmu)}$. This goal is reached through the employment of \emph{snapshots}, i.e. properly chosen solutions of \eqref{FE_ocp}. The low-dimensional function space is a subset of $\mathbb X \disc$ and in order to solve the optimality system for a new value of $\bmu \in \mathcal P$ it is sufficient to perform a standard Galerkin projection over the reduced space. In the following, we assume to have already built the reduced space\footnote{The description of the algorithm used to build the spaces is postponed in Section \ref{POD}} for the global variable $X(\bmu)$, say $\mathbb X_N \subset\mathbb X \disc$. We stress that the \emph{reduced dimension} $N$ verifies $N \ll \Cal N$. The building procedure is part of the \emph{offline phase} where not only the basis functions are computed, but also the $\bmu$-independent quantities are assembled and stored. \\ Once the offline phase is concluded, the \emph{reduced optimality system} reads:
given $\boldsymbol \mu \in \Cal P$, find
$X_N (\boldsymbol \mu) \eqdot (y_N(\boldsymbol \mu), u_N(\boldsymbol \mu),z_N(\boldsymbol \mu)) \in
\mathbb X_N \eqdot
\Cal \state_N \times \Cal \control_N \times \Cal \state_N$ such that it holds:
\begin{equation}
\label{ROM_ocp}
\begin{cases}
D_{y}\Lg(X_N; y_\text{d}, \bmu)[\omega] = 0 & \forall \omega \in \state{_N},\\
D_u\Lg(X_N; y_\text{d}, \bmu)[\kappa] = 0 & \forall \kappa \in \control{_N},\\
D_z\Lg(X_N; y_\text{d}, \bmu)[\zeta] = 0 & \forall \zeta \in \state{_N}.\\
\end{cases}
\end{equation}
The solution is given in an \emph{online phase}, where for every new evaluation of $\bmu \in \Cal P$, the optimality system \eqref{ROM_ocp} is assembled and solved. From the latest arguments, it is natural to deduce that one of the main ingredients of this procedure is the efficient division between a possibly expensive offline phase, which is performed only once, and a fast projection phase to solve the system for several parameters. This may be possible if we assume an affine decomposition for \eqref{ROM_ocp}, i.e. the involved equations have the following form:
\begin{equation}
\label{affinity}
\begin{aligned}
D_{y}\Lg(X_N, y_\text{d}, \bmu)[\omega] =
\displaystyle \sum_{\mathsf q=1}^{Q_y} \Theta_{{y}}^\mathsf q(\boldsymbol{\mu})D_{y}\Lg^\mathsf q(X_N; y_\text{d})[\omega],\\
D_{u}\Lg(X_N; y_\text{d}, \bmu)[\kappa] =
\displaystyle \sum_{\mathsf q=1}^{Q_u} \Theta_{{u}}^\mathsf q(\boldsymbol{\mu})D_{u}\Lg^\mathsf q(X_N; y_\text{d})[\kappa],\\
D_{z}\Lg(X_N; y_\text{d}, \bmu)[\zeta] =
\displaystyle \sum_{\mathsf q=1}^{Q_z} \Theta_{{z}}^\mathsf q(\boldsymbol{\mu})D_{z}\Lg^\mathsf q(X_N; y_\text{d})[\zeta].
\end{aligned}
\end{equation}
In other words, the system can be recast as the product of $\boldsymbol \mu -$dependent smooth functions  \\$\Theta_{{y}}^\mathsf q(\boldsymbol{\mu}), \Theta_{u}^\mathsf q(\boldsymbol{\mu}), \Theta_{z}^\mathsf q(\boldsymbol{\mu})$ and  $\boldsymbol{\mu} -$independent forms
$D_{y}\Lg^\mathsf q(X_N, y_\text{d})[\omega], D_{u}\Lg^\mathsf q(X_N; y_\text{d})[\kappa],$ and \\$D_{z}\Lg^\mathsf q(X_N; y_\text{d})[\zeta]$. \\
When this is the case, the online phase does not depend on $\Cal N$ and usually guarantees the solution of the system in a small amount of time. 
\begin{remark}
\B{For nonlinear systems, even if structure \eqref{affinity} is fulfilled, the  involved nonlinear forms still depend on $X^{\mathcal N}(\bmu)$ and this affects the computational advantage in using a POD-Galerkin approach, since it involves the assembly (and projection) of the high fidelity solution \Aa{during the online stage}. To overcome this issue, hyper-reduction techniques based on the Empirical Interpolation Method (EIM) may be employed, see e.g. \cite{barrault2004empirical} or \cite[Chapter 5]{hesthaven2015certified}.}
\end{remark}
In the next Section, we will describe the space-time POD-Galerkin strategy to build the reduced space $\mathbb X_N$.

\subsection{Offline and Online phase: from Space-Time POD Algorithm for OCP($\boldsymbol{\mu}$)s to Galerkin Projection}
\label{POD}
The ROM building process has mainly been addressed through {two} techniques: the POD \cite{ballarin2015supremizer, burkardt2006pod, Chapelle2013, hesthaven2015certified} and the greedy algorithm \cite{gerner2012certified, hesthaven2015certified, negri2015reduced, negri2013reduced, rozza2007stability}. In this work, we will focus on the first one since it can be applied to any state equation: indeed, greedy constructions are based on the employment of {an error estimator, which is} still not available for nonlinear time dependent \ocp s.
\\ The POD-Galerkin algorithm samples $N_{\text{max}}$ parameters in $\mathcal P$ and computes the related snapshots. After this \emph{exploratory phase}, a compressing stage starts, where $ N < N_{\text{max}}$ basis functions are provided after snapshots manipulation, aiming to get rid of the redundant information in the parametrized system. \\
Let us define the subset $\Cal P_{N_{\text{max}}} \subset \Cal P$ given by the sampled parameters, which has cardinality $N_{\text{max}}$. The chosen snapshots will form the following \emph{sampled} manifold
$$
\mathcal M_{N_{\text{max}}}^{\Cal N}=  \{ X \disc(\boldsymbol{\mu}) \;| \; \boldsymbol{\mu} \in \Cal P_{N_{\text{max}}}\} \subset \mathcal M^{\Cal N},
$$
where we assume $N_{\text{max}}$ large enough to let $\mathcal M_{N_{\text{max}}}^{\Cal N}$ be a reliable representation $\mathcal M^{\Cal N}$. The POD algorithm is applied separately for each involved variable, in a \emph{partitioned space-time approach}, where the procedure provides spaces of dimension $N$ \B{which minimizes the following quantities}:
\begin{equation*}
\label{crit}
\hspace{-1.5cm}
\sqrt{\frac{1}{N_{\text{max}}}
\sum_{\boldsymbol{\mu} \in \Cal P_{N_{\text{max}}}} \underset{\omega_N \in {{\mathbb Y}}_{N}}{\text{min }} \norm{y\discy(\boldsymbol{\mu}) - \omega_N}_{\mathbb Y}^2},
\hspace{.2cm}
\sqrt{\frac{1}{N_{\text{max}}}
\sum_{\boldsymbol{\mu} \in \Cal P_{N_{\text{max}}}} \underset{\kappa_N \in {{\mathbb U}}_{N}}{\text{min }} \norm{u\discu(\boldsymbol{\mu}) - \kappa_N}_{\mathbb U}^2},
\end{equation*}
\begin{equation*}
\sqrt{\frac{1}{N_{\text{max}}}
\sum_{\boldsymbol{\mu} \in \Cal P_{N_{\text{max}}}} \underset{\zeta_N \in {{\mathbb Y}}_{N}}{\text{min }} \norm{z\discy(\boldsymbol{\mu}) - \zeta_N}_{\mathbb Y}^2}.
\end{equation*}
In the following we describe the snapshots manipulation to build the reduce space only for one variable, say the state $y (\boldsymbol \mu)$. The proposed arguments can be identically replied for the other variables as well. First of all, we consider the ordered set of parameters $\boldsymbol{\mu}_1, \dots, \boldsymbol{\mu}_{N_{max}}\in \Cal P_{N_{\text{max}}}$ to which correspond a set of order snapshots  $y \discy(\boldsymbol{\mu}_1), \dots, y \discy(\boldsymbol{\mu}_{N_{max}})$. The correlation matrix $\mbf C^{y } \in \mbb R^{N_{max} \times N_{max}}$ of snapshots of the state variable, i.e.:
$$
\mbf C_{ml}^{y } = \frac{1}{N_{max}}({y }\B{\discy}(\boldsymbol{\mu}_m),y \discy(\boldsymbol{\mu}_l))_{\mathbb {Y}}, \hspace{1cm} 1 \leq m,l \leq N_{max}.
$$
First, we solved the following eigenvalue problem:
$$
\mbf C^{y} x_n^{y } = \lambda_n^{y } x_n^{y }, \hspace{1cm} 1 \leq n \leq N,
$$
with $\norm {x_n^{y }}_{\mathbb {Y}} = 1$.
Let us assume to have sorted the eigenvalues $\lambda_1^{y }, \dots, \lambda_{N_{\text{max}}}^{y }$  in decreasing order and to \Aa{consider} only the first $N$ ones, namely $\lambda_1^{y }, \dots, \lambda_{N}^{y }$, and the corresponding eigenvectors
$x_1^{y }, \dots, x_{N}^{y }$.
Let $(x_n^{y })_m$ be \emph{m-th} component of the state eigenvector $x_n^{y } \in \mbb R^{N_{max}}$. Thus, the \Aa{POD basis functions} are given by the following relation:
\begin{equation}
\label{eq:basis}
\chi_n^{y } = \displaystyle \frac{1}{\sqrt{{\lambda_{\Aa{n}}^{y }}}}\sum_{m = 1}^{\B{N_{max}}} (x_n^{y })_m y\discy (\boldsymbol{\mu}_m), \hspace{1cm} 1 \leq n \leq N.
\end{equation}
\B{The relation \eqref{eq:basis} is standard in data-compression algorithms such as POD, see e.g.  \cite{hesthaven2015certified, quarteroni2015reduced} }. \B{ It represents how the bases can be written in terms of the POD eigenvalues-eigenvectors pairs.}
We remark that the time instances are not separated in the POD procedure: i.e. the snapshots preserve the space-time structure. \B{The values $N_{max}$ and $N$ can be guided by the analysis of the POD eigenvalues, since the following holds \cite{hesthaven2015certified, quarteroni2015reduced}:
\begin{equation}
\sqrt{\frac{1}{N_{\text{max}}}
\sum_{m = 1}^{N_{\text{max}}}  \norm{y\discy(\boldsymbol{\mu}_{m}) - P_N(y\discy(\boldsymbol{\mu}_m)) }_{\mathbb Y}^2} = \sqrt{
\sum_{m = N + 1}^{N_{\text{max}}}\lambda_m^y,}
\end{equation}
where $P_N: \mathbb Y \rightarrow \mathbb Y_N$ is the projection for functions in $\mathbb Y$ onto the reduced space $\mathbb Y_N$.} \B{Another aspect to take care of is the sampling of the $N_{max}$ parameters for the POD, which can be related to some previous knowledge one has of the system at hand.}\\
As already specified in Section \ref{FE}, the linearized \ocp $\,$ leads to the solution of a saddle point system at each iteration of the Newton's method. To prove the existence and uniqueness of the solution, the matrix $\mathsf B$ of system \eqref{J_saddle} must verify the inf-sup stability condition \eqref{FE_lbb} for every $\boldsymbol \mu \in \Cal P$. 
The relation, at the space-time level, holds thanks to the same discretization technique used for state and adjoint variable. However, this assumption does not guarantees the fulfillment of the inf-sup stability at the reduced level. Indeed, the basis has to be manipulated in order to achieve this goal since the standard space-time POD process may lead to different reduced spaces for state and adjoint, even if the high fidelity discretization is the same for both the variables. To overcome this issue, we exploit \emph{aggregated spaces} technique. This strategy is very common and well known in ROM literature for \ocp s, see \cite{bader2016certified,bader2015certified,dede2010reduced,gerner2012certified,karcher2014certified,karcher2018certified, negri2015reduced,negri2013reduced,quarteroni2007reduced} as references. The main idea is to build a common function space which represents both state and adjoint variables
\begin{equation}
\label{state_r}{\state}_{N}= \text{span }\{\chi^{y}_n, \chi^{z}_n, \; n = 1, \dots, N\},
\end{equation}
while for the control variable we rely on the standard space
\begin{equation}
\label{control_r} {\control}_{N} = \text{span}\{\chi^{u}_n, \; n = 1, \dots, N\}.
\end{equation}
We can now define the \emph{basis} matrix
$$
\mathsf Z =
\begin{bmatrix}
\mathsf Z_{\mathsf x} \\
\mathsf Z_{\mathsf z}
\end{bmatrix},
\spazio \text{and} \spazio
\mathsf Z_{\mathsf x} =
\begin{bmatrix}
\mathsf Z_{\mathsf y} \\
\mathsf Z_{\mathsf u}
\end{bmatrix}
$$
where
$
\mathsf Z_{\mathsf y} \equiv \mathsf Z_{\mathsf z} = [\chi_{1}^{y} | \cdots | \chi_{N}^{y}| \chi_{1}^{z} | \cdots | \chi_{N}^{z}] \in \mathbb R^{ \Cal N_{y} \times 2N}
$ and $
\mathsf Z_{\mathsf u} = [\chi_{1}^{u} | \cdots | \chi_{N}^{u}] \in \mathbb R^{\Cal N_{u} \times N}.
$ The $\mathsf Z $ spans the reduced space $\mathbb X_N$.
In this framework we can solve the optimality system in a small amount of time for every parametric instance through a Galerkin projection into the reduced spaces. Thus, the final system reads:
\begin{equation}
\label{G_compact_ROM}
\mathsf G_{N}(\mathsf X_N; \boldsymbol \mu) \mathsf X_N = \mathsf F_N,
\end{equation}
with
$$\mathsf G_{N}(\mathsf X_N; \boldsymbol \mu) \eqdot \mathsf Z^T \mathsf G(\mathsf Z \mathsf X_N; \boldsymbol \mu),
\spazio \text{and} \spazio \mathsf F_N \eqdot \mathsf Z^T \mathsf F.$$
The system \eqref{G_compact_ROM} inherits the nonlinearity feature from the high fidelity one and the Newton's method can be employed also in this case
\begin{equation}
\mathsf {X}_N^{j + 1} \eqdot \mathsf {X}_N^j+ \mathsf{Jac}_N(\mathsf X_N^{j}; \boldsymbol \mu)^{-1}(\mathsf F_N - \mathsf G(\mathsf X_N^j; \boldsymbol \mu)\mathsf X_N^j), \spazio j \in \mathbb N,
\end{equation}
presenting the saddle point structure for the Frech\'et derivative, i.e.
\begin{equation}
\label{Frechet_ROM}
\mathsf{Jac}_N (\mathsf X_N; \boldsymbol \mu) \mathsf X_N =
\begin{bmatrix}
\mathsf A_N & \mathsf B_N^T \\
\mathsf B_N & 0 \\
\end{bmatrix}
\begin{bmatrix}
\mathsf x_N \\
\mathsf z_N
\end{bmatrix},
\end{equation}
with $\mathsf{Jac}_N (\mathsf X_N; \boldsymbol \mu) = \mathsf Z^T \mathsf{Jac}(\mathsf Z \mathsf X_N; \bmu)\mathsf Z$, $\mathsf A_N = \mathsf Z_{\mathsf x}^T \mathsf A \mathsf Z_{\mathsf x}$ and $\mathsf B_N = \mathsf Z_{\mathsf z}^T \mathsf B \mathsf Z_{\mathsf x}$ \\
In the reduced framework, the \emph{reduced inf-sup condition} reads
\begin{equation}
\label{ROM_infsup}
\beta_{N}(\bmu) \eqdot \adjustlimits \inf_{0 \neq \mathsf z_N} \sup_{0 \neq \mathsf x_N} \frac{\mathsf z_N^T\mathsf B_N \mathsf x_N}{\norm{\mathsf x_N}_{\mathbb Y \times \mathbb U}\norm{\mathsf z_N}_{\Cal \state}} \geq \hat{\beta}_{N}(\bmu) > 0,
\end{equation}
which is verified due to the aggregated spaces definition. \Bb{It is known that using the state and adjoint function spaces built on their respective snapshots without any kind of manipulation might be not sufficient to guarantee the well-posedness of a linear saddle point system \eqref{Frechet_ROM}. Indeed, by means of standard POD function spaces, the fulfillment of \eqref{ROM_infsup} is not assured. To avoid this issue, we follow the strategy already employed in \cite{negri2013reduced} and previously in \cite{dede2010reduced}: we use the same space for state and adjoint, made by the union of the two different bases obtained by the POD applied to state and adjoint, respectively, as defined in \eqref{state_r}. We remark that, even if the reduced system increases its final dimension through this approach, it is usually much smaller then $\Cal N$.} \Aaa{Moreover, it is clear that to have a good reduced approximation, the high fidelity solution must be a good representation of the continuous framework. However, as already specified in Section \ref{3}, it has to be verified case by case. This kind of analysis goes beyond the goal of this contribution and we will always assume that the space-time approximation is a valid discrete representation of the continuous problem}. In the next Section we will show with some numerical examples on how convenient ROMs can be in the framework of space-time nonlinear \ocp s.

\section{Application to Shallow Waters Equations}
\label{5}
This Section shows how the previously proposed methodologies can be applied to the \A{viscous} SWEs model. We stress that the described strategies are general and can be applied to several state equations. Furthermore, they easily adapt to simpler settings like steady and/or linear problems, which are still useful in many engineering and scientific fields, as already specified in Section \ref{intro}. However, in the following, we will focus on environmental sciences and, more specifically, in coastal management: a growing impact field of research that needs the support of accurate real-time simulations. To this end, first, we will introduce the \A{viscous} SWEs state equations and some motivations for the primary importance of such a model in coastal fluid dynamics. Then, we will briefly describe the optimality system formulation at the continuous level before proposing some numerical results to test our methodology in a physical and geometrical parametrized setting.

\subsection{Main Motivations and Problem Formulation}
\label{SWE}
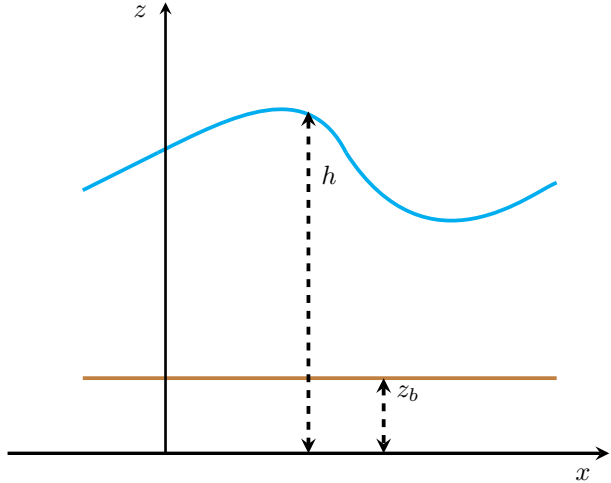
\begin{figure}[H]
\begin{center}
\begin{tikzpicture}
\draw[line width=0.5mm,  cyan] (0,0) .. controls (1,0.5) and (2,1) .. (2.5,0);
\draw [line width=0.5mm,  cyan](2.5,0) .. controls (3.6,-1.7) and (5,-0.5) .. (5.3,-0.4);
\draw[line width=0.5mm,  cyan] (0,0) .. controls (-1,-0.5) and (-1,-0.5) .. (-1,-0.5);

\draw[line width=0.5mm,  brown] (-1,-3)  -- (5.3,-3);

\draw [arrow, line width=0.35mm] (0.1,-4) -- (0.1, 2);
\draw [arrow, line width=0.45mm] (-2,-4) -- (6, -4);
\draw[b_arrow_dashed,  line width=0.5mm] (2,0.55) -- (2, -4);
\draw[b_arrow_dashed,  line width=0.5mm] (3,-3.) -- (3, -4 );
\node[anchor=south, label=left:{$\vspace{2mm}x$}] (6,-4.4) at (6,-4.4) {};
\node[anchor=north, label=left:{$\vspace{2mm}z$}] (0.1,2) at (0.1,2) {};
\node[anchor=west, label=left:{$\vspace{2mm}h$}] (2.5,- 0.3) at (2.5,- 0.3) {};
\node[anchor=west, label=left:{$\vspace{2mm}z_b$}]  (3.6,-3.5) at (3.6,-3.2) {};
\end{tikzpicture}
\end{center}
\caption{Notations: schematic representation.}
\label{notations}
\end{figure}

This Section is motivated by the growing demand of fast and reliable simulations in the framework of coastal management. Indeed, the marine environment is related to social and economic growth, biodiversity and ecosystem preservation, monitoring plans for possibly dangerous events related to weather factors or anthropic behaviour.
\\ The \A{viscous} SWEs result in a very versatile model in coastal sciences which can represent marine impact on shores and coasts, planetary currents, tsunamis waves... \cite{cavallini2012quasi,Shallow}. The state equation and the control problem have been studied analytically and numerically in several works, see for example \cite{agoshkov1994recent2,  Agoshkov1993, agoshkov1994recent, agoshkov2007optimal, ferrari2004new, miglio1999finite, miglio2005model, miglio2005model2, ricchiuto2007application, ricchiuto2009stabilized,  takase2010space}. In a parametrized setting, they have been explored as state equation in \cite{Navon2, Navon1} and in \cite{Strazzullo3} in an optimal control framework. \\ First of all we define our parameter, $\bmu = (\mu_1, \mu_2, \mu_3, \mu_4) \in \mathcal P \subset \mathbb R^4$. The first three parameters are related to the physics of the problem, while the latter describes the geometry of the considered spatial domain: the specific description of the parametrized setting is postponed later in the Section. For this peculiar state equation, we define the function spaces $Y_v = H^1_{\Gamma _{D_v}(\mu_4)}(\Omega(\mu_4))$, $Y_h = L^2(\Omega(\mu_4))$ and the space  $U = L^2(\Omega(\mu_4))$, where $\Gamma_{D_v}(\mu_4)$ is a portion of the boundary domain $\partial \Omega(\mu_4)$ where Dirichlet boundary conditions have been imposed. The involved variables are the vertically averaged velocity profile of the wave $\boldsymbol v$ and the surface elevation variable $h$, respectively considered in $Y_v$ and $Y_h$. We will focus on the simpler setting, where the bottom bathymetry $z_b$ is defined as a constant function: it can be generalized to more realistic bathymetry, see e.g \cite{ferrari2004new, ricchiuto2009stabilized, ricchiuto2007application}. In Figure  \ref{notations} we provide a description of the physical phenomenon we are dealing with, together with the notations we will exploit.  We used the standard 2D-model presented in \cite{Agoshkov1993, miglio1999finite}. 
The state solution valocity and height are defined, respectively, in
$$\mathbb Y_v = \Big \{ \boldsymbol v \in L^2(0,T; [Y_v]^2) \text{ such that } \dt{\boldsymbol v} \in L^2(0,T; [Y_v\dual]^2) \Big \},$$ and
$$\mathbb Y_h = \Big \{ h \in L^2(0,T; Y_h) \text{ such that } \dt{h} \in L^2(0,T; Y_h\dual) \Big \}.$$
Thus, we can define the global state function space given by $\mathbb Y = \mathbb Y_v \times \mathbb Y_h$, where the state variable $(\boldsymbol v, h)$ will be sought. Furthermore, let  $\boldsymbol u \in \mathbb U \eqdot L^2(0,T; [U]^2)$ be the control variable of the system. For the numerical results, we will deal with a distributed optimal control, where the control variable acts as an external forcing term such as atmospheric pressure, bottom friction and wind action.  Namely, you are not actually controlling the system but one can interpret the optimal control framework as an inverse problem capable to guess what are the physical conditions that can represent a desired velocity-height profile $(\boldsymbol v_d, h_d) \in \mathbb Y_{\text{obs}} \eqdot \mathbb Q_v \times \mathbb Q_h$, where 
$\mathbb Q_v  \eqdot L^2(0,T; [L^2(\Omega(\mu_4))]^2)$ and $\mathbb Q_h  \eqdot L^2(0,T; L^2(\Omega(\mu_4)))$. The problem at hand reads: given $\bmu \in \mathcal P$ find $(\boldsymbol v, h)\in \mathbb Y$ which minimizes  $  J((\boldsymbol v, h), \boldsymbol u, (\boldsymbol v_{\text{d}}, h_{\text{d}}))$, where 
\begin{align}
\label{J}
J((\boldsymbol v, h), \boldsymbol u; (\boldsymbol v_{\text{d}}, h_{\text{d}})) = &  
\half \int_{0}^T\int_{\Omega(\mu_4)}{(h - h_{\text d}(\mu_3))^2d\Omega(\mu_4)dt} \\ \nonumber
& \quad \quad \quad + \half \int_{0}^T\int_{\Omega(\mu_4)}{(\boldsymbol v - \boldsymbol v_{\text d}(\mu_3))^2d\Omega(\mu_4)dt} +
\alf \int_{0}^T\int_{\Omega(\mu_4)}{\boldsymbol u^2d\Omega(\mu_4)dt},
\end{align}
constrained to the following equation 
\begin{equation}
\label{SWEs_state}
\begin{cases}
\displaystyle \dt{\boldsymbol v} - \mu_1 \Delta \boldsymbol v + \mu_2(\boldsymbol v \cdot \nabla)\boldsymbol v + g \nabla h - \boldsymbol u = 0 & \text{in } \Omega(\mu_4) \times (0, T) , \vspace{1mm}\\
\displaystyle \dt{h} + \dive(h \boldsymbol v) = 0 & \text{in } \Omega(\mu_4) \times (0,T), \\
\boldsymbol v = \boldsymbol v_0 & \text{on }  \Omega(\mu_4) \times \{  0 \}, \\
h = h_0 & \text{on }   \Omega(\mu_4) \times \{  0 \}, \\
\boldsymbol v = \boldsymbol 0& \text{on }  \partial  \Omega(\mu_4) \times (0, T).
\end{cases}
\end{equation}
\B{We recall that $\boldsymbol v \in L^2(0,T; [Y_v]^2)$, $\boldsymbol v_t \in L^2(0,T; [{Y_v}^*]^2)$, $h \in L^2(0,T; Y_h)$, $h_t \in L^2(0,T; {Y_h}^*)$ and $\boldsymbol u  \in L^2(0,T; [U]^2)$, where $Y_h = U = L^2(\Omega(\mu_4))$ and $Y_v = H^1_{\Gamma_D}(\Omega(\mu_4))$. Now, noticing that
$$
Y_v \hookrightarrow U \hookrightarrow Y_v^{\dual},
$$
we can interpret the first equation as a sum of terms in $L^2(0,T; [{Y_v}^{*}]^2)$. We now focus on the second equation. Using standard regularity results for Navier-Stokes equations \cite{constantin1988navier}, $h$ is in $L^2(0,T; Y_v)$. Exploting Sobolev embedding theorems in dimension \Aa{three (considering space and time)}, we have $h \boldsymbol v \in L^2(0,T; [Y_h]^2)$ and, consequently, $\dive{(h\boldsymbol v)} \in L^2(0,T; Y_v^*)$. Thus, the second equation is a sum of terms in
$L^2(0,T; {Y_v}^{*})$, since $Y_h\dual \hookrightarrow Y_v^{\dual}$.}
\B{Even if the setting presented in Section \ref{3} includes a broad class of state equations, to the best of our knowledge no results about \ocp s governed by the two-dimensional SWEs are known in the Lagrangian context. However, we believe that one could apply the same techniques already used in \cite[Section 1.8.1]{hinze2008optimization} for Navier-Stokes equations to state the well-posedness of the problem and to recover hypothesis (i) - (ix). However, this topic goes beyond the goal of this contribution and we will restrict ourselves to the case of well-posedness for this specific test case.}
As already specified in Section \ref{3}, the value of $\alpha \in (0,1]$ influences the control action: the smaller is $\alpha$, the larger the action of the control variable is.
The state equation \eqref{SWEs_state} models free surface incompressible flows in  hydrostatic pressure: this assumption is verified when shallow depths are considered, namely when the water height is much lower than the wavelength, which is classical for coastal phenomena. We deal with two physical parameters $\mu_1$ and $\mu_2$ representing diffusive and convective action of the system considered, respectively. Moreover, $\mu_3$ will define the desired solution profile, while $\mu_4$, as already specified, will affect the spatial geometry we will deal with. 
Let us suppose to have recast the state equation \eqref{SWEs_state} in weak formulation as 
$\mathcal E((\boldsymbol v, h), \boldsymbol u; \bmu)$, thus we can define
 the following Lagrangian functional 
\begin{equation}
\label{Lagrangian}
\Lg((\boldsymbol v, h), \boldsymbol u,  (\boldsymbol z, q), (\boldsymbol v_{\text d}, h_{\text d}); \bmu) = J((\boldsymbol v, h), \boldsymbol u; (\boldsymbol v_{\text{d}}, h_{\text{d}})) + \la (\boldsymbol z, q), \mathcal E((\boldsymbol v, h), \boldsymbol u; \bmu)\ra_{\mathbb Q \mathbb Q\dual},
\end{equation}
where the variable $ (\boldsymbol z, q) \in \mathbb Y$ is the adjoint variable. Here we do not specify all the bilinear forms involved, the interested reader may refer to \cite{Strazzullo3}. The optimality system is, then, obtained through the differentiation with respect to the all involved variables. The minimization problem thus reads: given $\bmu \in \mathcal P$, find 
$((\boldsymbol v, h), \boldsymbol u,  (\boldsymbol z, q)) \in \mathbb Y \times \mathbb U \times \mathbb Y$ such that
\begin{equation}
\label{optimality_system_SWEs}
\begin{cases}
D_{\boldsymbol v}\Lg((\boldsymbol v, h), \boldsymbol u, (\boldsymbol z, q))[\boldsymbol \zeta] = 0 & \forall \boldsymbol \zeta \in \mathbb Q_v,\\
D_h\Lg((\boldsymbol v, h), \boldsymbol u, (\boldsymbol z, q))[\theta] = 0 & \forall \theta \in \mathbb Q_h,\\
D_{\boldsymbol u}\Lg((\boldsymbol v, h), \boldsymbol u, (\boldsymbol z, q))[\boldsymbol \tau] = 0 & \forall \boldsymbol \tau \in \mathbb U,\\
D_{\boldsymbol z} \Lg((\boldsymbol v, h), \boldsymbol u, (\boldsymbol z, q)) [\boldsymbol \kappa]= 0 & \forall \boldsymbol \kappa \in \mathbb Q_v,\\
D_q\Lg((\boldsymbol v, h), \boldsymbol u, (\boldsymbol z, q))[\xi] = 0 & \forall \xi \in \mathbb Q_h.\\
\end{cases}
\end{equation}
For the sake of completeness, we report the system \eqref{optimality_system_SWEs} in strong formulation:
\begin{align}
\label{SWEs_opt}
\begin{cases}
\displaystyle \boldsymbol v - \dt{\boldsymbol z} - \mu_1 \Delta \boldsymbol z - \mu_2 (\boldsymbol v \cdot \nabla)\boldsymbol z + \mu_2 (\nabla \boldsymbol v)^T\boldsymbol z - h\nabla q = \boldsymbol v_d(\mu_3) & \text{in } \Omega(\mu_4) \times (0, T), \vspace{1mm}\\
\displaystyle  h - \dt{q} - \boldsymbol v \cdot \nabla q - g\dive(\boldsymbol z) = h_d(\mu_3) & \text{in } \Omega(\mu_4) \times (0, T), \vspace{1mm}\\
\boldsymbol z = \boldsymbol 0  & \text{on } \partial \Omega(\mu_4) \times (0, T) \\
\boldsymbol z = \boldsymbol 0 & \text{on } \Omega(\mu_4) \times \{  T \}, \\
q= 0 & \text{on }  \Omega(\mu_4) \times \{  T \}, \\
\alpha \boldsymbol u = \boldsymbol z & \text{in } \Omega(\mu_4) \times (0, T) \\
\displaystyle \dt{\boldsymbol v} - \mu_1 \Delta \boldsymbol v + \mu_2 (\boldsymbol v \cdot \nabla)\boldsymbol v + g \nabla h = \boldsymbol u & \text{in } \Omega(\mu_4) \times (0, T), \vspace{1mm}\\
\displaystyle \dt{h} + \dive(h \boldsymbol v) = 0 & \text{in } \Omega(\mu_4) \times (0, T), \\
\boldsymbol v = \boldsymbol 0  & \text{on } \partial \Omega(\mu_4) \times (0, T),\\
\boldsymbol v = \boldsymbol v_0 & \text{on } \Omega(\mu_4) \times \{  0 \}, \\
h = h_0 & \text{on }  \Omega(\mu_4) \times \{  0 \}. \\
\end{cases}
\end{align}
In the next Section we will present the numerical results for a distributed \ocp s with physical and geometrical parameters. They are an extension of \cite{Strazzullo3} to parametrized spatial domain. Furthermore, we will briefly discuss the high fidelity space-time approximation and how we perform a POD-Galerkin in this specific framework in order to build the reduced space approximation for this specific \ocp s.
\subsection{Numerical Results}

\label{Results}
\begin{figure}[H]
\centering

\includegraphics[width=0.45\textwidth]{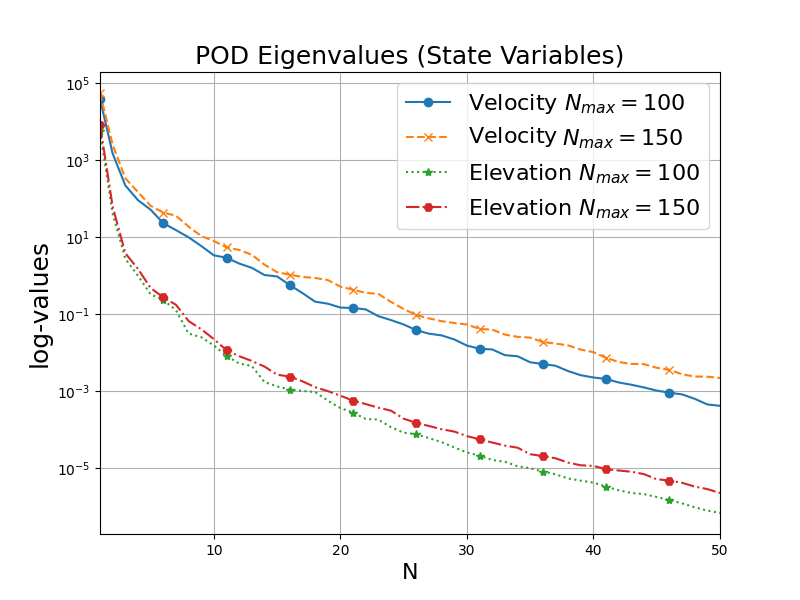}
\includegraphics[width=0.45\textwidth]{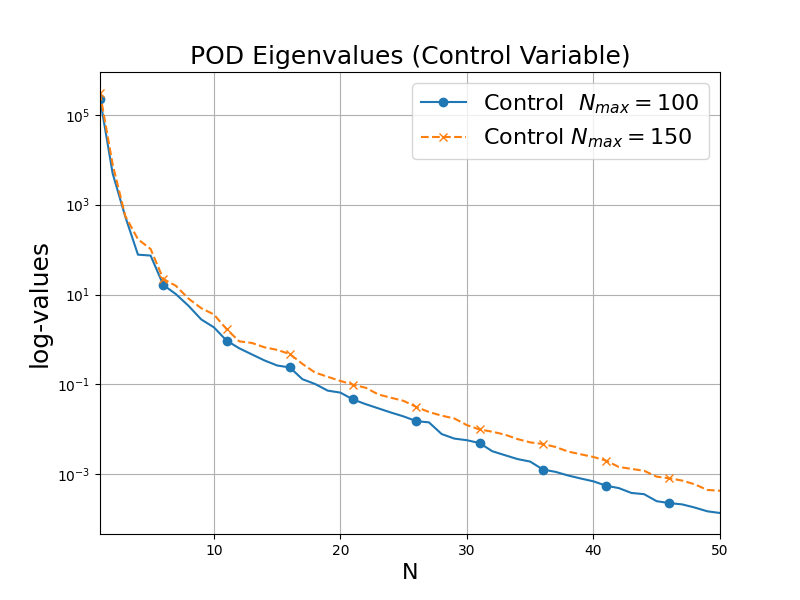}\\
\includegraphics[width=0.45\textwidth]{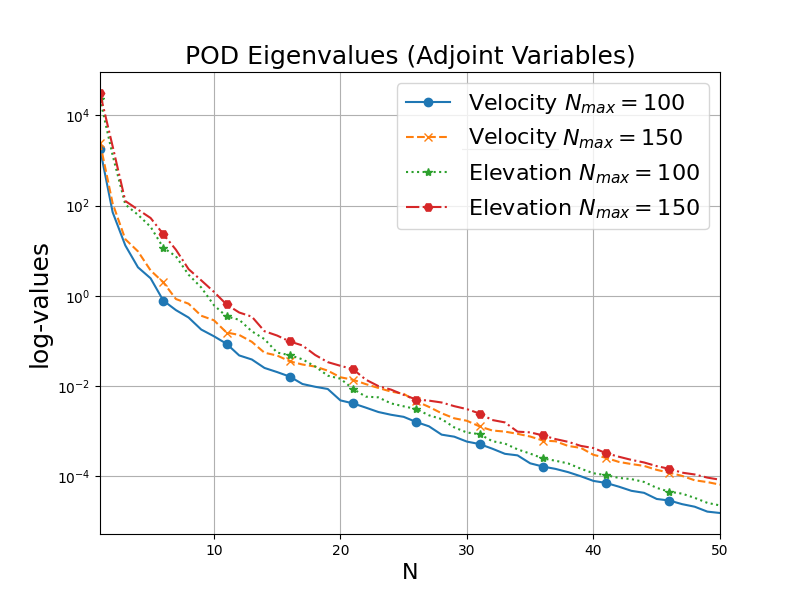}
\caption{Comparison of POD eigenvalues for $N_{\text{max}}=100$ and $N_{\text{max}}=150$. \emph{Top Left}: state variables. \emph{Top Right}: control variable. \emph{Bottom}: adjoint variables.  }

\label{fig:eigs}
\end{figure}

This Section aims at validating the performances of the POD-Galerkin projection reduced approch proposed in Section \ref{4}. We followed the test case presented in \cite{ferrari2004new}, recasting it in the parametrized setting briefly introduced in Section \ref{SWE}, where not only physical parameter is considered, but also the gemetry of the spatial domain changes with respect to the choise of $\bmu \in  \Cal P = (0.00001,1.) \times (0.01, 0.5) \times (0.1, 1.)  \times (0.8, 1.5)$. Indeed, the spatial domain is given by $\Omega(\mu_4) = [0, 10\mu_4] \times [0, 10]$. In an inverse problem fashion, the optimal control setting gives information about the forcing term needed to achieve a given desired solution profile. To perform simulations and consequent reduction, we took into account a pull back of the optimality system \eqref{optimality_system_SWEs} in the \emph{reference domain} given by $\Omega = [0, 10] \times [0, 10]$, i.e. the reference parameter is $\mu_4 = 1$.
In the considered framework, we exploited $z_b = 0$ for the flat bathymetry and we studied the wave evolution in the time interval $[0,T]$ where the final time is $T = 0.8s$. Let $x_1$ and $x_2$ be the spatial coordinates. The problem at hand aims at reducing the impact of the spreading of a mass of water with an initial Gaussian distributed elevation under a null initial velocity: i.e.
$$
\boldsymbol v_0 = \boldsymbol 0, \spazio \text{and} \spazio h_0 = 0.2(1 + 5e^{(-(\frac{x_1}{\mu_4} - 5)^2 - (x_2 - 5)^2 + 1))}).
$$
Namely we want to study under which physical conditions, say wind action and bottom friction, the system can reach a the desired state $(\mu_3 \boldsymbol v_{\text d}, \mu_3 h_{\text d})$, where  $(\boldsymbol v_{\text d}, h_{\text d})$ is the solution at time $T$ of the uncontrolled state equation \eqref{SWEs_state}, with
$$\boldsymbol {v_{\text d}}_0 = \boldsymbol 0, \spazio \text{and} \spazio {h_{\text d}}_{0} = 2e^{(-(\frac{x_1}{\mu_4} - 5)^2 - (x_2 - 5)^2 + 1)},$$
and no forcing term, i.e. $\boldsymbol u = \boldsymbol 0$. \\
Let us briefly analyse the high fidelity approximation. We exploited FE discretization, and the optimal solution fields have been discretized through linear polynomial, namely $N^y_{\text{FE}} = N^v_{\text{FE}} + N^h_{\text{FE}}$, where  $ N^v_{\text{FE}}$ and  $ N^h_{\text{FE}}$ are the FE dimensions for state and adjoint velocity and height profile and control variable, respectively, obtained using $r_v = r_h = 1$, as proposed in \cite{saleri2007geometric}. The same polynomial degree has been exploited for the control variable, giving a FE discretization of dimension $N^u_{\text{FE}}$. Concerning the time approximation, we perform Euler methods, dividing the time interval $N_t = 8$ timesteps, with a $\Delta t = 0.1s$. The number of timesteps can be increased following the iterative techniques presented in \cite{HinzeStokes, HinzeNS, Stoll1, Stoll}. Although, in this work, for the sake of simplicity and clearness, we exploited a direct solver for the optimality system \eqref{SWEs_opt} and this affects the resolution of the high fidelity approximation. In the end, at the truth approximation level, we deal with a system of a total dimension ${\Cal {N}} = 76352$, although, we underline that it is sufficient to validate the reduced approach proposed in Section \ref{4}. Let us focus our attention on the construction phase of the low dimensional function spaces. As specified in Section \ref{POD}, we applied a partitioned approach, where several POD compressions, over the five different correlation matrices of $N_{\text{max}} = 100$ snapshots have been carried out for all the involved variables.\B{ The number of snapshots $N_{\text{max}}$ was heuristically chosen observing that the POD eigenvalue decay is comparable to the one obtained by $N_{\text{max}}=150$, see Figure \ref{fig:eigs}. Besides, we could not explore further the parametric space since the space-time approximation drastically affects the computational time needed for the offline phase of the POD algorithm, i.e.\hspace{-1mm} the snapshots collection.}
Then, the application of the strategy described in Section \ref{POD} will lead to the following spaces:
\begin{align*}
& {\mathbb Y}^{v}_{N} = \text{span}\{\chi^{v}_n, \; n = 1, \dots, N\}, \\
& {\mathbb Y}^{h}_{N} = \text{span}\{\chi^{h}_n, \; n = 1, \dots, N\}, \hspace{2mm}\\
& {\mathbb U}_{N} = \text{span}\{\chi^{u}_n, \; n = 1, \dots, N\}, \\
& {\mathbb Y}^{z}_{N} = \text{span}\{\chi^{ z}_n, \; n = 1, \dots, N\}, \\
& {\mathbb Y}^{ q}_{N} = \text{span}\{\chi^{q}_n, \; n = 1, \dots, N\}, \\
\end{align*}
where we retained $N = 30$ basis functions\footnote{For the sake of simplicity, we exploit the same value of $N$ for all the spaces.}, given a global reduced space dimension of $9N = 270$. \B{The value $N$ has been heuristically chosen as a trade-off between relative errors values and computational time saving in the reduced framework. }Thus, we define
$$
{\mathbb Z}^{vz}_{N}= \text{span }\{\chi^{v}_n, \chi^{z}_n, \; n = 1, \dots, N\}
\hspace{2mm} \text{and} \hspace{2mm}
{\Cal Z}^{hq}_{N}= \text{span }\{\chi^{h}_n, \chi^{q}_n, \; n = 1, \dots, N\}.
$$
For state and adjoint velocity-height variable we will use $\mathbb Z_N = {\mathbb Z}^{vz}_{N} \times {\mathbb Z}^{hq}_{N}$ so that the inf-sup condition \eqref{ROM_infsup} holds.  Indeed, we remark that, this strategy guarantees the reduced inf-sup stability condition \eqref{ROM_infsup} to be verified for the linearized problem at each step of the Newton's Method.
Despite the increasing of the dimensionality, we will see that the performances of the reduced projection is still convenient with respect to the high fidelity simulations: comments on the computational advantages are postponed later in the Section. \A{ We stress that the used POD-Galerkin approach might be sub-optimal for the state equation at hand as underlined in several works, see e.g. \cite{GRIMBERG, TaddeiSWE}}. \A{However, the proposed optimal control numerical setting does not suffer much the reduced representation: indeed we are dealing with a viscous model and, furthermore, the controlled framework seems to act itself as a sort of stabilization, where the forcing term changes in order to achieve a less convection-dominated solution.}
\Bb{This is confirmed both by the comparison between projection-based ROM and best-fit projection averaged relative log-errors for a testing set of 20 uniformly distributed parameters of Table \ref{tb:table_proj} and by the eigenvalue decay represented in Figure \ref{fig:eigs}.

\begin{table}

\caption{\Bb{Projection-based ROM and best-fit projection averaged relative log-errors comparison.}}
\label{tb:table_proj} 
\begin{center}\footnotesize{
\begin{tabular}{|c||c|c|c|c|c|c|c|c|c|c|}
\hline
 & \multicolumn{5}{c|}{\textbf{Averaged relative errors proj-ROM}} & \multicolumn{5}{c|}{\textbf{Averaged relative errors best-fit projection}} \\  \cline{2-11}
\multirow{-2}{*}{$N$} &$\boldsymbol v$&$h$ &$\boldsymbol u$ & $\boldsymbol w$ & $q$ & $\boldsymbol v$&$h$ &$\boldsymbol u$ & $\boldsymbol w$ & $q$ \\ 
 \hline
$6$ & 4.68e--2 &5.42e--3 & 2.30e--2 &5.26e--2 & 2.39e--2 & 2.52e--2 & 4.05e--3 & 1.17e--2 & 2.11e--2 & 2.21e--2 \\ \hline
$10$ &  1.48e--2 & 2.28e--3 & 4.53e--3 & 1.28e--2 & 7.34e--3 & 1.22e--2 & 1.58e--3 & 4.47e--3 & 1.04e--2 & 5.89e--3 \\ \hline
$16$ & 6.82e--3 & 9.01e--4 & 1.97e--3 & 5.72e--3 & 3.27e--3 & 5.65e--3 & 5.88e--4 & 1.97e--3 & 4.73e--3 &  2.50e--3 \\ \hline
$20$ & 4.14e--3 & 5.83e--4 & 1.34e--3 & 3.96e--3 & 2.07e--3 & 3.44e--3 & 3.89e--4 & 1.34e--3 & 3.40e--3 & 1.52e--3 \\ \hline
$26$ & 2.43e--3 & 3.17e--4 & 7.94e--4 & 4.47e--3 & 1.16e--3 & 2.07e--3 & 2.04e--4 & 7.88e--4 & 2.13e--3 & 8.39e--4 \\ \hline
$30$ & 1.81e--3 & 2.25e--4& 5.48e--4 & 1.75e--3 & 7.73e--4 & 1.15e--3 & 1.44e--4 & 5.48e--4 & 1.52e--3 & 5.54e--4 \\ \hline
\end{tabular}}
\end{center}
\end{table}}

\begin{figure}[H]
\centering

\includegraphics[width=0.45\textwidth]{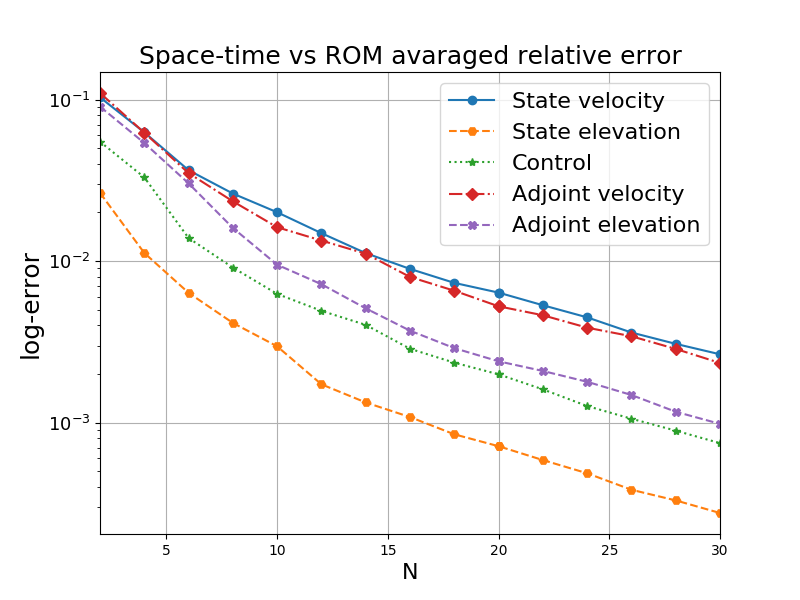}
\caption{Averaged relative error between the space-time and ROM approximation for all the involved variables.}

\label{err}
\end{figure}

All the features of the offline and online phase of the experiment are reported in Table \ref{tab:test_1}.
\begin{table}[H]
\centering
\caption{Data for the distributed \ocp  $\,$ governed by SWEs.}
\label{tab:test_1}
\footnotesize{
\begin{tabular}{|c|c|}
\hline
 \textbf{Data} & \textbf{Values} \\
\hline
$\mathcal P$ &  $(0.00001,1) \times (0.01, 0.5) \times (0.1, 1) \times (0.8, 1.5)$\\
\hline
$[0,T]$ &  $[0s, 0.8s]$\\
\hline
 values of $(\mu_1, \mu_2, \mu_3,\Aa{\mu_4},  \alpha)$ &  $\B{ (0.1, .01, .1, 1.5, 0.1)}$ \\
\hline
$N_{\text{max}}$ & 100\\
\hline
$N$ & 30 \\
\hline
Sampling Distribution & Uniform \\
\hline
$\Cal N$ & $76352$ \\
\hline
ROM System Dimension & 270 \\
\hline

\end{tabular}
}
\end{table}
The basis function considered allowed us to well describe the full order approximated system in the reduced framework, as the reader can notice from the average relative errors represented in Figure \ref{err} with the following norms:
$$
\intTime{\norm{\boldsymbol v ^{N_{\text{FE}}^v}- \boldsymbol v_N}^2_{H^1}}, \hspace{2mm}
\hspace{2mm}
\intTime{\norm{h^{N_{\text{FE}}^h} - h_N}^2_{L^2}}, \hspace{2mm}
$$
$$
\intTime{\norm{\boldsymbol u^{N_{\text{FE}}^u} - \boldsymbol u_N}^2_{L^2}}, \hspace{2mm}
\hspace{2mm}
\intTime{\norm{\boldsymbol z ^{N_{\text{FE}}^v} - \boldsymbol z_N}^2_{H^1}}, \hspace{2mm}
\text{and }\hspace{2mm}
\intTime{\norm{q^{N_{\text{FE}}^h} - q_N}^2_{L^2}}.\hspace{2mm}
$$

\begin{figure}[H]
\centering
\begin{subfigure}[b]{0.3\textwidth}
\centering
\includegraphics[width=0.8\textwidth]{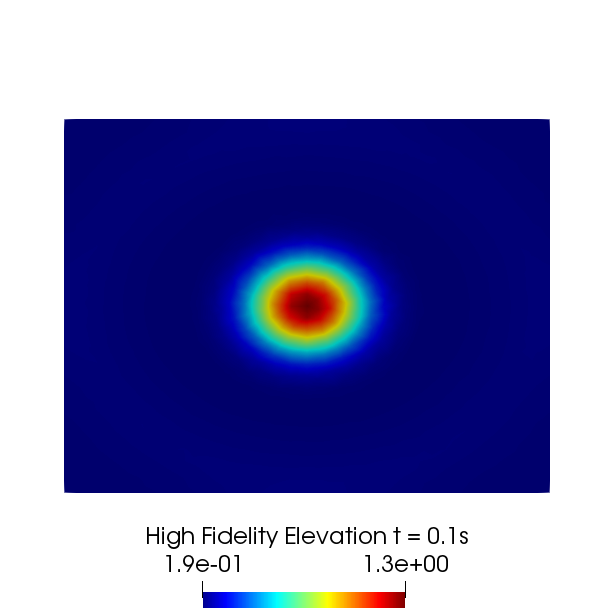}
\caption{}
\label{fig:off_h_1}
\end{subfigure}
\hfill
\begin{subfigure}[b]{0.3\textwidth}
\centering
\vspace{2mm}
\includegraphics[width=0.8\textwidth]{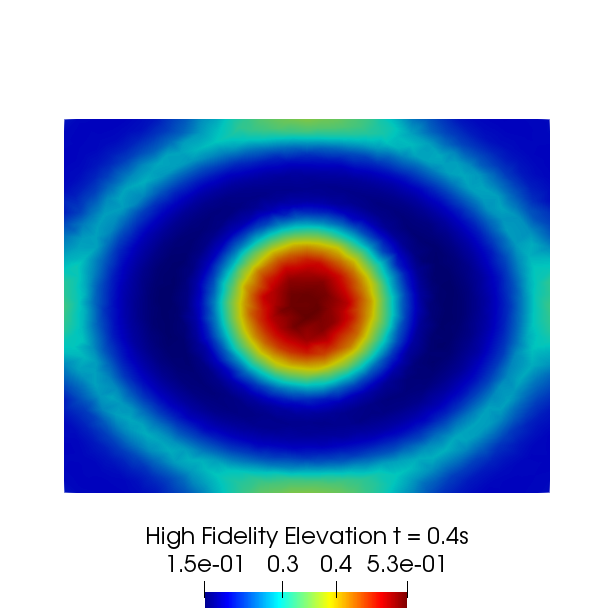}
\caption{}
\label{fig:off_h_4}
\end{subfigure}
\hfill
\begin{subfigure}[b]{0.3\textwidth}
\centering
\includegraphics[width=0.8\textwidth]{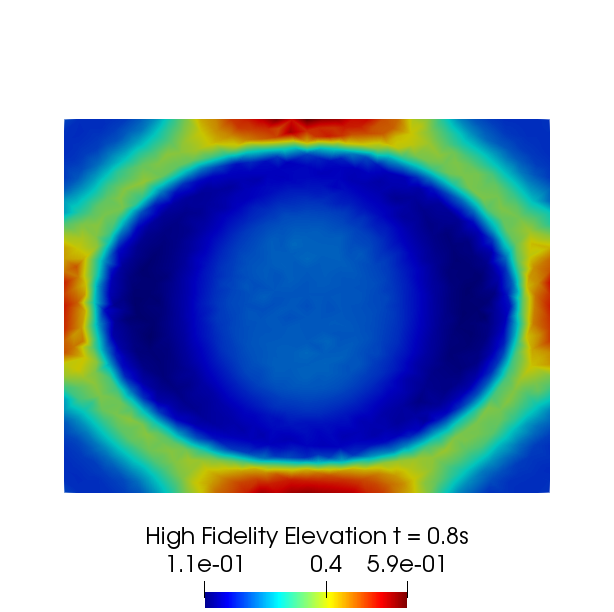}
\caption{}
\label{fig:off_h_8}
\end{subfigure}
\hfill
\begin{subfigure}[b]{0.3\textwidth}
\centering
\includegraphics[width=0.8\textwidth]{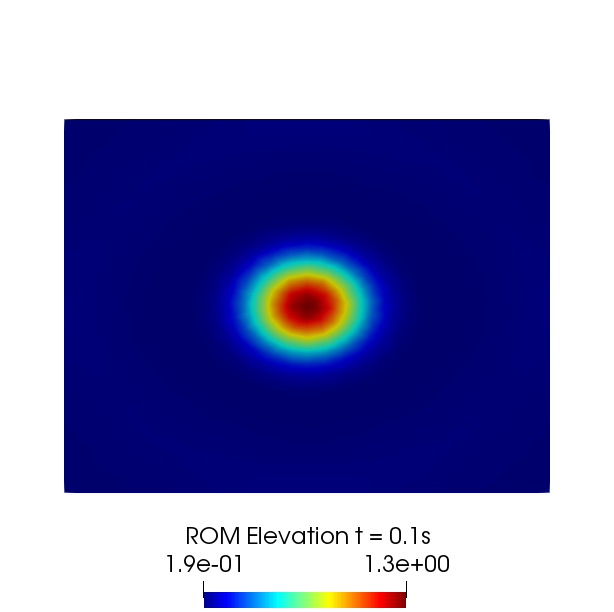}
\caption{}
\label{fig:on_h_1}
\end{subfigure}
\hfill
\begin{subfigure}[b]{0.3\textwidth}
\centering
\vspace{2mm}
\includegraphics[width=0.8\textwidth]{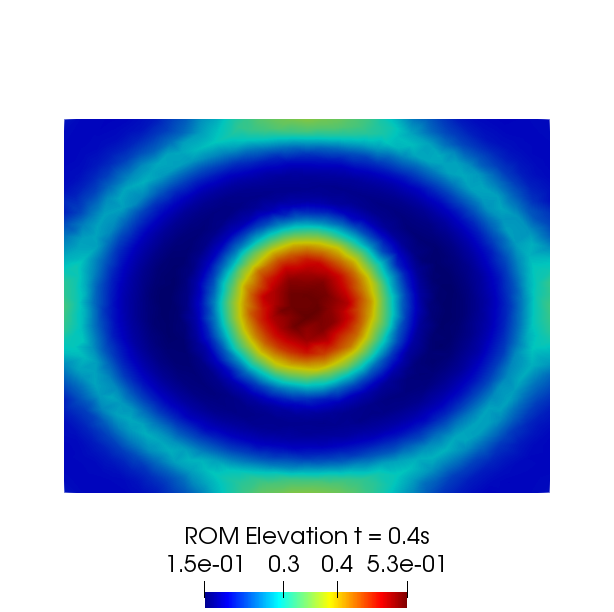}
\caption{}
\label{fig:on_h_4}
\end{subfigure}
\hfill
\begin{subfigure}[b]{0.3\textwidth}
\centering
\includegraphics[width=0.8\textwidth]{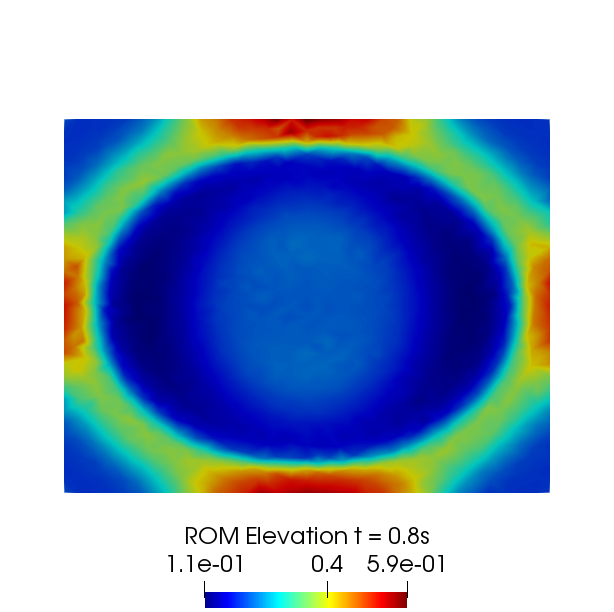}
\caption{}
\label{fig:on_h_8}
\end{subfigure}
\caption{Optimal high fidelity and reduced state elevation variable with $\alpha = 0.1$ and $\bmu= (0.1, .01, .1, 1.5)$. High fidelity solutions for $t = 0.1s, 0.4s, 0.8s$ in (a), (b), (c), respectively,  and reduced solutions for $t = 0.1s, 0.4s, 0.8s$ in (d), (e), (f).}
\label{h_comp}
\end{figure}

\begin{figure}[H]
\centering
\begin{subfigure}[b]{0.3\textwidth}
\centering
\includegraphics[width=0.8\textwidth]{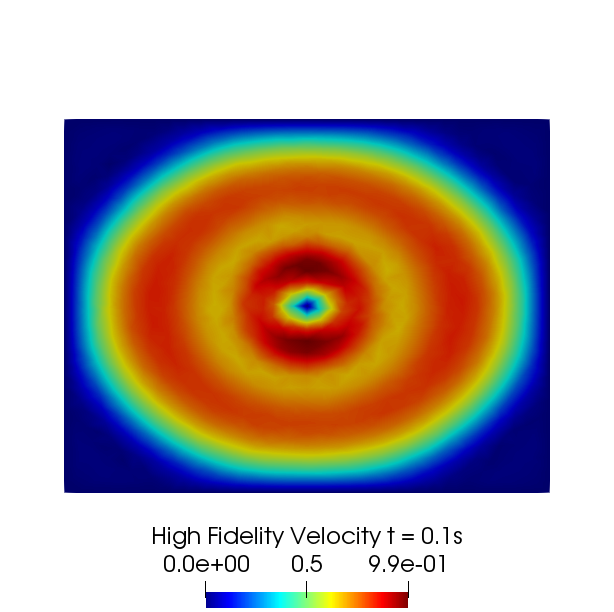}
\caption{}
\label{fig:off_v_1}
\end{subfigure}
\hfill
\begin{subfigure}[b]{0.3\textwidth}
\centering
\vspace{2mm}
\includegraphics[width=0.8\textwidth]{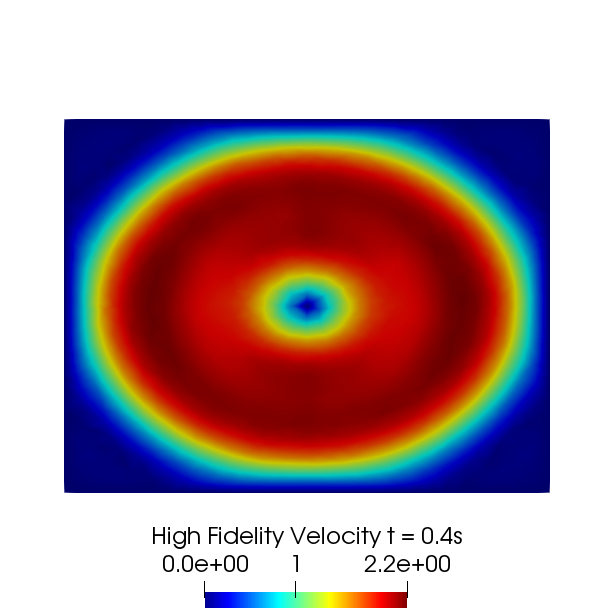}
\caption{}
\label{fig:off_v_4}
\end{subfigure}
\hfill
\begin{subfigure}[b]{0.3\textwidth}
\centering
\includegraphics[width=0.8\textwidth]{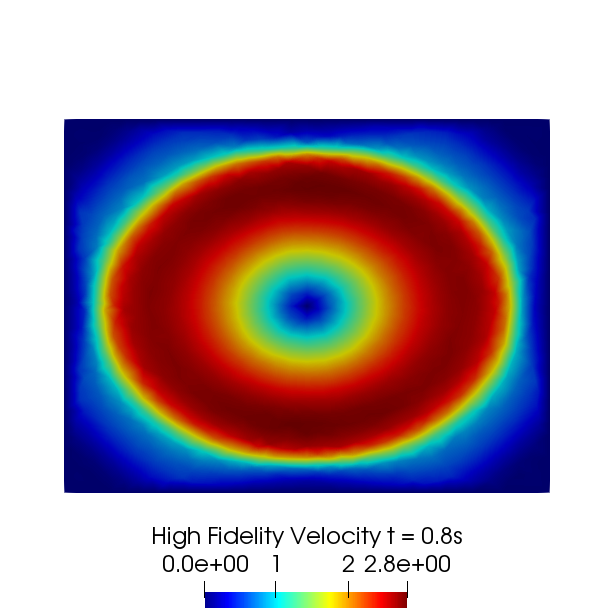}
\caption{}
\label{fig:off_v_8}
\end{subfigure}
\hfill
\begin{subfigure}[b]{0.3\textwidth}
\centering
\includegraphics[width=0.8\textwidth]{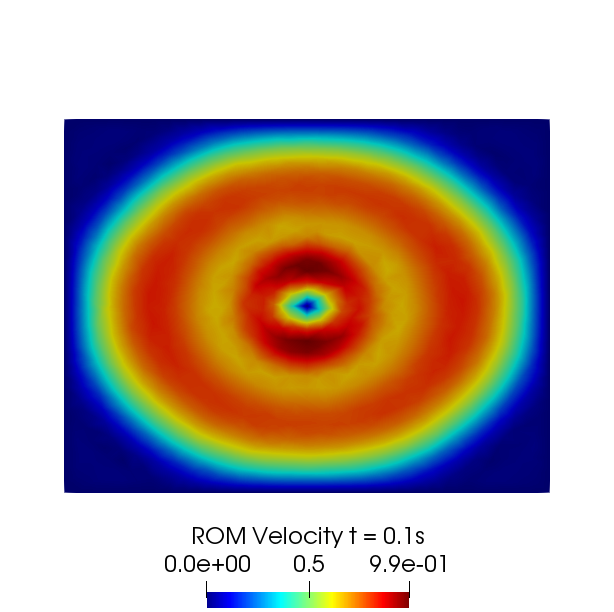}
\caption{}
\label{fig:on_v_1}
\end{subfigure}
\hfill
\begin{subfigure}[b]{0.3\textwidth}
\centering
\vspace{2mm}
\includegraphics[width=0.8\textwidth]{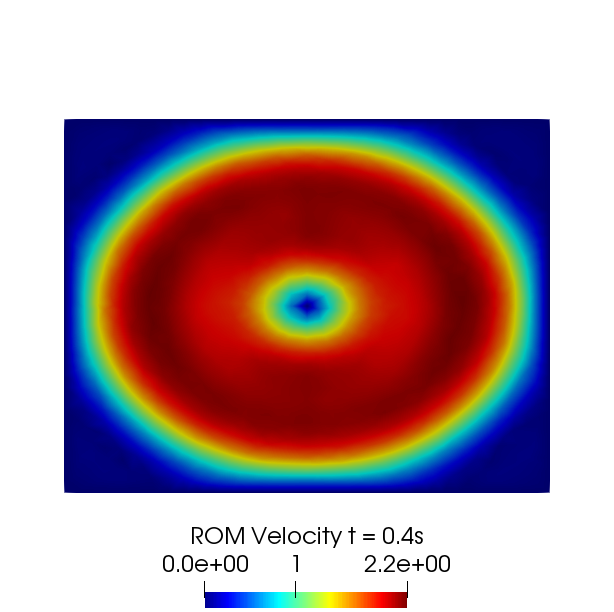}
\caption{}
\label{fig:on_v_4}
\end{subfigure}
\hfill
\begin{subfigure}[b]{0.3\textwidth}
\centering
\includegraphics[width=0.8\textwidth]{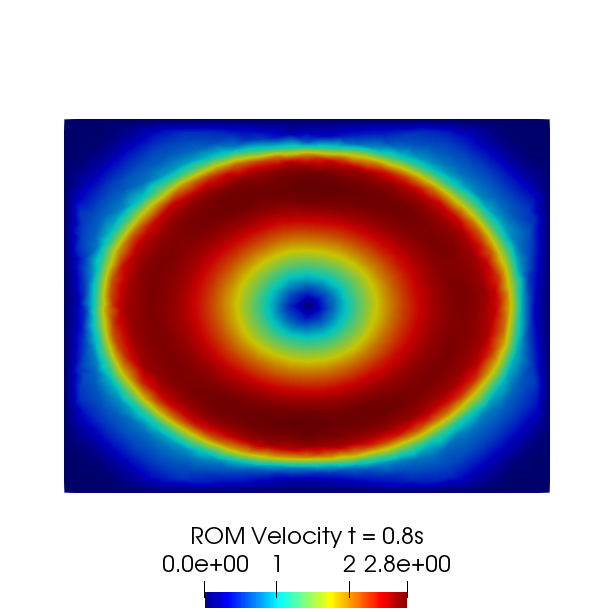}
\caption{}
\label{fig:on_v_8}
\end{subfigure}
\caption{Optimal high fidelity and reduced state velocity variable with $\alpha = 0.1$ and $\bmu= (0.1, .01, .1, 1.5)$. High fidelity solutions for $t = 0.1s, 0.4s, 0.8s$ in (a), (b), (c), respectively,  and reduced solutions for $t = 0.1s, 0.4s, 0.8s$ in (d), (e), (f).}
\label{v_comp}
\end{figure}

\begin{figure}[H]
\centering
\begin{subfigure}[b]{0.3\textwidth}
\centering
\includegraphics[width=0.8\textwidth]{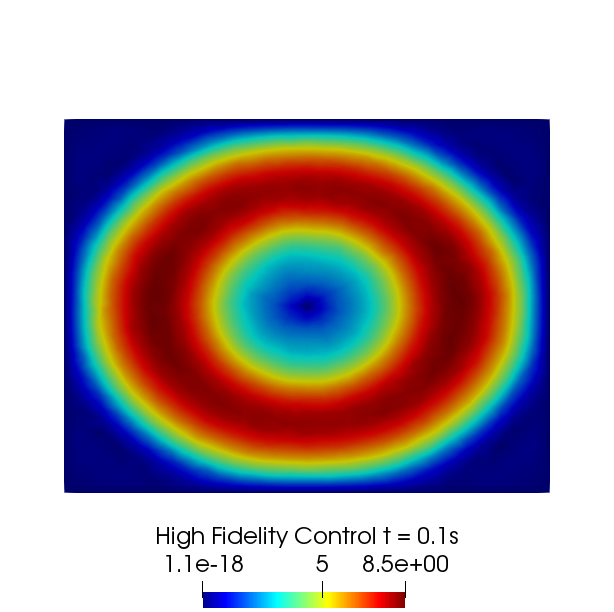}
\caption{}
\label{fig:off_u_1}
\end{subfigure}
\hfill
\begin{subfigure}[b]{0.3\textwidth}
\centering
\vspace{2mm}
\includegraphics[width=0.8\textwidth]{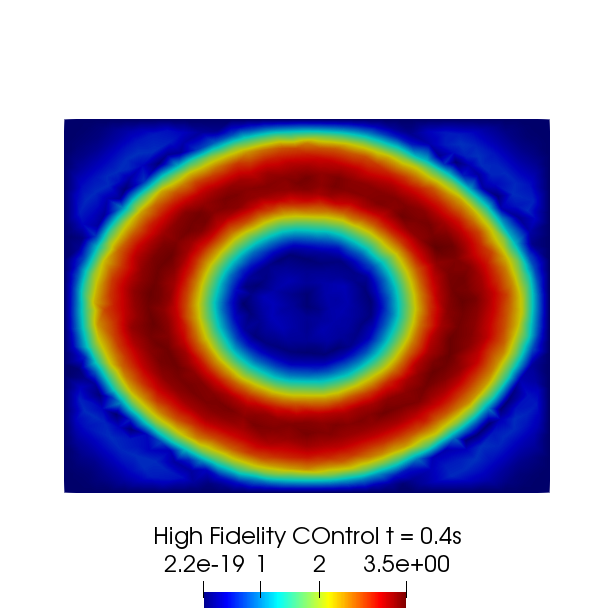}
\caption{}
\label{fig:off_u_4}
\end{subfigure}
\hfill
\begin{subfigure}[b]{0.3\textwidth}
\centering
\includegraphics[width=0.8\textwidth]{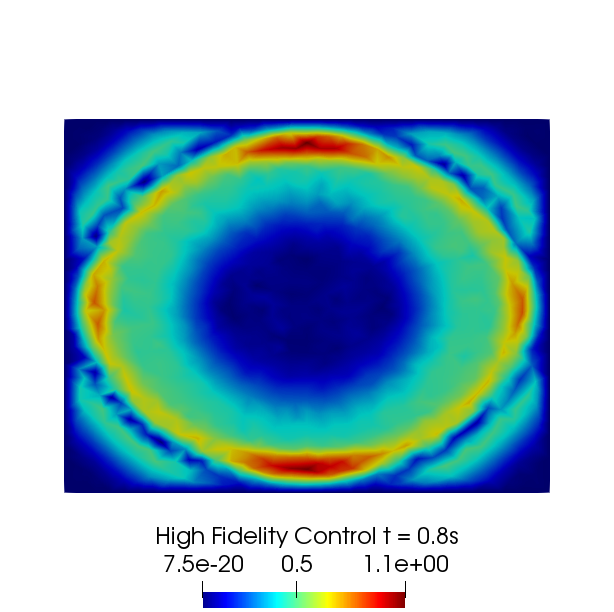}
\caption{}
\label{fig:off_u_8}
\end{subfigure}
\hfill
\begin{subfigure}[b]{0.3\textwidth}
\centering
\includegraphics[width=0.8\textwidth]{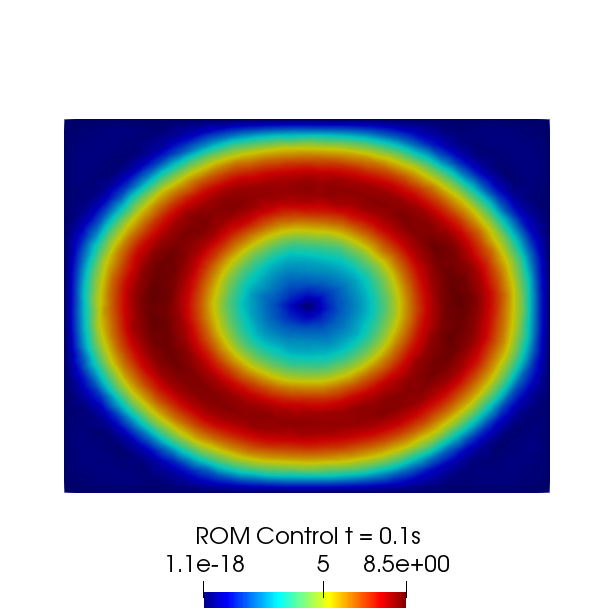}
\caption{}
\label{fig:on_u_1}
\end{subfigure}
\hfill
\begin{subfigure}[b]{0.3\textwidth}
\centering
\vspace{2mm}
\includegraphics[width=0.8\textwidth]{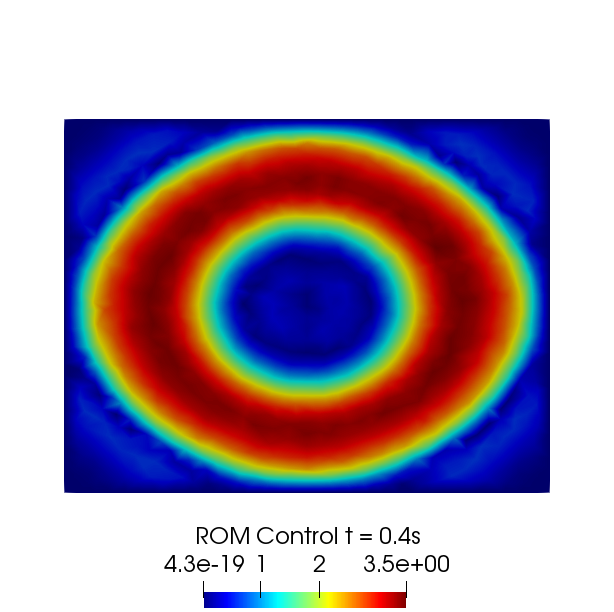}
\caption{}
\label{fig:on_u_4}
\end{subfigure}
\hfill
\begin{subfigure}[b]{0.3\textwidth}
\centering
\includegraphics[width=0.8\textwidth]{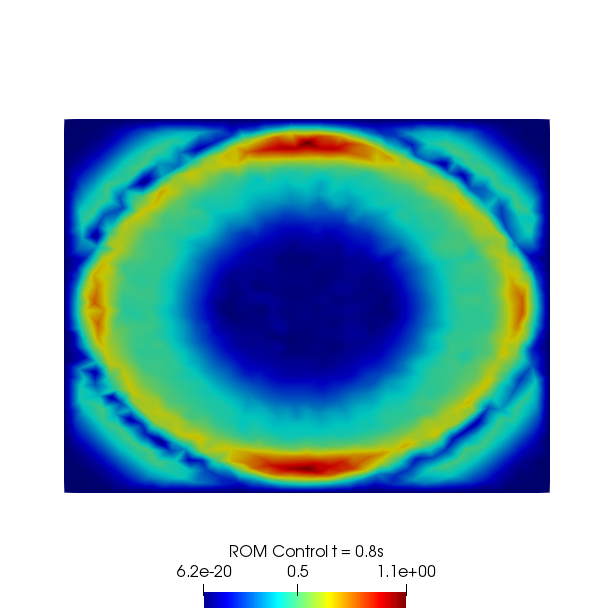}
\caption{}
\label{fig:on_u_8}
\end{subfigure}
\caption{Optimal high fidelity and reduced control variable with $\alpha = 0.1$ and $\bmu= (0.1, .01, .1, 1.5)$. High fidelity solutions for $t = 0.1s, 0.4s, 0.8s$ in (a), (b), (c), respectively,  and reduced solutions for $t = 0.1s, 0.4s, 0.8s$ in (d), (e), (f).}
\label{u_comp}
\end{figure}

The errors are averaged over a testing set of $20$ parameters uniformly distributed: from the plots in Figure \ref{err}, one can observe how the POD-Galerkin approach leads to a good approximation of all the involved quantities. Indeed, using $N = 30$, state and adjoint velocity reach values around $2\cdot 10^{-3}$, while the relative error for control, state and adjoint elevation is below $10^{-3}$.  Furthermore, the accuracy of the reduced model can be understood also from the comparison between the space-time solutions and the ROM solutions presented in Figures \ref{v_comp},
\ref{h_comp} and \ref{u_comp} for state velocity, state elevation and control at $t = 0.1s, 0.4s, 0.8s$, respectively.
Let us comment on the computational time performances between space-time approximation and ROM simulations. Indeed, not only the reduced basis are able to reproduce several time instances in an accurate way, but there is a gain in the computational time needed for parametrized simulations. 

We call \emph{speedup index} the number of reduced simulations which can be performed in the time of the solution of one high fidelity optimality system. For this specific test case, the speedup is of the order of $O(20)$ for $N = 1, \dots, 30$, namely, it is lightly influenced by the reduced spaces dimension.
The speedup index tells us that performing a Galerkin projection in the aggragated spaces is still convenient with respect the solution of the whole space-time OCP($\bmu$). The saved computational time can be used to study and analyse several parametric configurations in a real-time or many-query context. \\
In the next Section some comments and conclusions follow.

\section{Conclusions}

\label{6}
In this contribution, we propose ROMs as a reliable and fast strategy to deal with parametrized nonlinear time dependent \ocp $\:$. The analysis we propose is general and can be applied to several PDE($\bmu$)s constrained optimization processes. We described the optimality system used to reach a desired profile both at the continuous and discrete level, in a space-time fashion, presenting it not only from a theoretical point of view, but also in its algebraic form, underlining the hidden saddle point structure of the linearized system. Thus, we moved towards model order reduction in order to alleviate the issue of the complexity of the optimality system, which results in a high dimensional problem. We propose a space-time POD-Galerkin approach: the choice was led by the need of an algorithm which could be applied also to very complicated equations, such as the nonlinear time dependent optimality systems. This work relies in building a reduced framework for time dependent nonlinear \ocp s, capable of filling the gap between data and physical model. We validate the method in environmental applications such as marine ecosystem management and coastal engineering, though a solution tracking optimization problem governed by \A{viscous} SWEs under physical and geometrical parametrization. Our aim was to show how ROMs could be a suitable tool to rapidly simulate marine environment, deeply characterized by a growing demanding computational effort. Indeed, the general proposed methodology results in fast simulations without paying in accuracy with respect to the time consuming space-time approximation. \\


\section*{Acknowledgements}
We acknowledge the support by European Union Funding for Research and Innovation -- Horizon 2020 Program -- in the framework of European Research Council Executive Agency: Consolidator Grant H2020 ERC CoG 2015 AROMA-CFD project 681447 ``Advanced Reduced Order Methods with Applications in Computational Fluid Dynamics''. We also acknowledge the PRIN 2017  ``Numerical Analysis for Full and Reduced Order Methods for the efficient and accurate solution of complex systems governed by Partial Differential Equations'' (NA-FROM-PDEs) and the INDAM-GNCS project ``Tecniche Numeriche Avanzate per Applicazioni Industriali''.
The computations in this work have been performed with RBniCS \cite{rbnics} library, developed at SISSA mathLab, which is an implementation in FEniCS \cite{fenics} of several reduced order modelling techniques; we acknowledge developers and contributors to both libraries.

\bibliographystyle{abbrv}
\bibliography{maria}

\end{document}